\documentclass{amsart}
\newtheorem{theorem}{Theorem}[section]

\theoremstyle{definition}

\theoremstyle{remark}

\numberwithin{equation}{section}

\newcommand{\scr}{\scriptstyle}
\def\sumprime_#1{\setbox0=\hbox{$\scriptstyle{#1}$}
\setbox2=\hbox{$\displaystyle{\sum}$}
\setbox4=\hbox{${}'\mathsurround=0pt$}
\dimen0=.5\wd0 \advance\dimen0 by-.5\wd2
\ifdim\dimen0>0pt
\ifdim\dimen0>\wd4 \kern\wd4 \else\kern\dimen0\fi\fi
\mathop{{\sum}'}_{\kern-\wd4 #1}}
\font\ger=eufm10
\newcommand{\gs}{\hbox{\ger S}}
\def\ndiv{\not \hskip .03in \mid}

\begin{document}
\title[HIGHER CORRELATIONS OF DIVISOR SUMS III]{Higher Correlations of Divisor Sums Related to Primes III: $k$-Correlations}
\author{D. A. Goldston}
\address{Department of Mathematics and Computer Science, San Jose
State University, San Jose, CA 95192, USA}
\email{goldston@mathcs.sjsu.edu}
\thanks{The first author was supported in part by an NSF Grant}
\author{C. Y. Y{\i}ld{\i}r{\i}m}
\address{ Department of Mathematics, Bo\~{g}azi\c{c}i University, Istanbul 80815, Turkey}
\email{yalciny@boun.edu.tr}

\date{\today}


\begin{abstract}
We obtain the general $k$-correlations for a short divisor sum related to primes.
\end{abstract}

\maketitle
\section{  Introduction}

This is the third in a series of papers concerned with the calculation and applications of higher correlations of short divisor sums that approximate the von Mangoldt function $\Lambda(n)$. We return to the short divisor sum approximation of the first paper defined for $R\ge 1$ by
\begin{equation} \Lambda_R(n)= \sum_{\stackrel{\scr d |n}{\scr d\le R}}\mu(d) \log (R/d). \label{1.1} \end{equation}
The correlations we are interested in evaluating are
\begin{equation} \mathcal{ S}_k(N,\text{\boldmath$j$}, \text{\boldmath $a$}) =\sum_{n=1}^N \Lambda_R(n+j_1)^{a_1}\Lambda_R(n+j_2)^{a_2}\cdots \Lambda_R(n+j_r)^{a_r}\label{1.2}\end{equation}
and
\begin{equation}\tilde{\mathcal{S}}_k(N,\text{\boldmath$j$}, \text{\boldmath$a$}) =\sum_{n=1}^N \Lambda_R (n+j_1)^{a_1}\Lambda_R (n+j_2)^{a_2}\cdots \Lambda_R(n+j_{r-1})^{a_{r-1}}\Lambda(n+j_r) \label{1.3}\end{equation}
where $\text{\boldmath$j$} = (j_1,j_2, \ldots , j_r)$ and $\text{\boldmath$a$} = (a_1,a_2, \ldots a_r)$, the $j_i$'s are distinct integers, $a_i\geq 1$ and $\sum_{i=1}^r a_i = k$. In \eqref{1.3} we assume that $r\ge 2$ and take $a_r=1$. We also define
\begin{equation} \tilde{\mathcal{ S}}_1(N, \text{\boldmath$j$}, \text{\boldmath$a$}) =\sum_{n=1}^N\Lambda(n+j_1) \sim N \label{1.4}\end{equation}
 if $|j_1|\le N$ by the prime number theorem. In the first paper in this series we evaluated $\mathcal{S}_k(N,\text{\boldmath$j$},\text{\boldmath$a$})$ and $\tilde{\mathcal{S}}_k(N,\text{\boldmath$j$},\text{\boldmath$a$})$ for $1\le k \le 3$, and also showed how for any $k$ the two correlations are related. In this paper we will consider the general case.
For $ \text{\boldmath$j$} = (j_1,j_2,\ldots , j_r)$, where the $j_i$'s are distinct integers, we define the singular series
\begin{equation} \gs(\text{\boldmath$j$}) = \prod_p\left( 1- \frac{1}{p}\right)^{-r}\left(1-\frac{\nu_p(\text{\boldmath$j$})} {p}\right) \label{1.5}\end{equation}
where $\nu_p(\text{\boldmath$j$})$ is the number of distinct residue classes modulo $p$ that the $j_i$'s occupy.

\begin{theorem} \label{Theorem1} Given $k\ge 1$,  let $\text{\boldmath$j$} = (j_1,j_2, \ldots , j_r)$ and $\text{\boldmath$a$} = (a_1,a_2, \ldots a_r)$, where the $j_i$'s are distinct integers, and $a_i\geq 1$ with $\sum_{i=1}^r a_i = k$. Assume $\max_i |j_i|\le R $ and $R\ge 2$. Then
we have
\begin{equation} \mathcal{S}_k(N,\text{\boldmath$j$},\text{\boldmath$a$}) = \big(\mathcal{ C}_k(\text{\boldmath$a$})\gs(\text{\boldmath$j$})+o_k(1)\big)N(\log R)^{k-r} +O(R^k),\label{1.6}\end{equation}
where the $\mathcal{ C}_k(\text{\boldmath$a$})$ are constants that are computable rational numbers. For $k\le 4$ we have
\begin{eqnarray}  \mathcal{ C}_1(1)&=&1, \nonumber \\
\mathcal{ C}_2(2)&=&1, \quad  \mathcal{ C}_2(1,1)=1, \nonumber \\
 \mathcal{ C}_3(3)&=& \frac{3}{ 4}, \quad \mathcal{ C}_3(2,1)=1, \quad \mathcal{ C}_3(1,1,1)=1,\nonumber \\
\mathcal{ C}_4(4)&=& \frac{3}{4}, \quad \mathcal{ C}_4(3,1)=\frac{3}{ 4}, \quad \mathcal{ C}_4(2,2)=1, \quad \mathcal{ C}_4(2,1,1) =1, \quad \mathcal{ C}_4(1,1,1,1) =1.\nonumber
\end{eqnarray}
Denoting $\mathcal{C}_k(k)$ as $\mathcal{C}_k$, all of the constants $\mathcal{ C}_k(\text{\boldmath$a$})$ with $\text{\boldmath$a$}=(a_1,a_2, \ldots , a_r)$ are determined from the constants $\mathcal{C}_k$ by the formula
\begin{equation} \mathcal{ C}_k(\text{\boldmath$a$})
= \prod_{i=1}^r \mathcal{C}_{a_i} .\label{1.7} \end{equation}
\end{theorem}
We actually prove a more precise result than \eqref{1.6}, but for applications the result stated here is sufficient.

If we wish to apply Theorem \ref{Theorem1} and obtain explicit numerical results we need to know the values of the  constants $\mathcal{C}_k$. These constants are defined in terms of  multiple integrals which may be evaluated by  residue calculations. These residue calculations becomes increasingly difficult as $k$ gets larger. We will obtain in this paper the values
\[ \mathcal{C}_1 = 1, \quad \mathcal{C}_2 = 1, \quad \mathcal{C}_3 = \frac{3}{4}, \quad \mathcal{C}_4 = \frac{3}{4} .\]
David Farmer has computed further values. He obtained the value $\mathcal{C}_4$ given above and also has found
\[  \mathcal{C}_5 = \frac{11065}{2^{14}} = .67535\ldots, \quad \mathcal{C}_6 = \frac{11460578803}{2^{34}}= .66709\ldots . \]

From the first paper in this series, we can evaluate the mixed correlations \eqref{1.3} as an immediate consequence of Theorem \ref{Theorem1}. We need the following form of the Bombieri-Vinogradov theorem. Let
\begin{equation} E(x;q,a)= \psi(x;q,a) - [(a,q)=1]\frac{x}{ \phi(q)},\label{1.8}\end{equation}
where $[P(x)]$ is the Iverson notation where square brackets around a true-false statement takes the value $1$ if the statement is true and $0$ if the statement is false. Suppose for some fixed $0<\vartheta \le 1$ that
\begin{equation} \sum_{1\le q\le x^{\vartheta - \epsilon}} \max_{\stackrel{\scr a }{ \scr (a,q)=1}}|E(x;q,a)|  \ll \frac{x}{\log ^\mathcal{A}x} \label{1.9} \end{equation}
holds for any $\epsilon >0$, any $\mathcal{A}=\mathcal{A}(\epsilon)>0$, and $x$ sufficiently large.  This is a weakened form of the Bombieri-Vinogradov theorem if $\vartheta=\frac{1}{2}$, and therefore \eqref{1.9} holds unconditionally if $\vartheta \le \frac{1}{2}$.  Elliott and Halberstam conjectured \eqref{1.9} is true with $\vartheta=1$. (The constant  $\vartheta$ is often referred to as the level of distribution of primes in arithmetic progression.)

\begin{theorem} \label{Theorem2} Given $ k\ge 2$,  let $\text{\boldmath$j$} = (j_1,j_2, \ldots , j_r)$ and $\text{\boldmath$a$} = (a_1,a_2, \ldots a_r)$, where $r\ge 2$, $a_r=1$, and where the $j_i$'s are distinct integers, $a_i\geq 1$, and $\sum_{i=1}^r a_i = k$. Assume $\max_i |j_i|\ll R^\epsilon$ and $R\gg N^\epsilon$. Then
we have, for $N^{\epsilon}\ll R \ll N^{\min(\frac{\vartheta}{k-1}, \frac{1}{k})-\epsilon}$ where \eqref{1.9} holds with $\vartheta$,
\begin{equation} \tilde{\mathcal{S}}_k(N,\text{\boldmath$j$},\text{\boldmath$a$}) = \big( \mathcal{C}_k(\text{\boldmath$a$})\gs(\text{\boldmath$j$})+o_k(1)\big) N(\log R)^{k-r} .\label{1.10}\end{equation}
where $\mathcal{ C}_k(\text{\boldmath$a$})$
are the same constants as in Theorem \ref{Theorem1}.\end{theorem}

We can also improve the range of $R$ in Theorem \ref{Theorem2}
to  $R \ll N^{\frac{\vartheta}{k-1} -\epsilon} $ by applying the methods of this paper directly to the mixed correlations, but this improvement is only significant in the conditional case when $ 1 - \frac{1}{k} \le \vartheta \le 1$. We will delay our applications to primes to the next paper in this series. To optimize some applications we will need to generalize slightly the previous results by allowing for different truncation lengths in the correlations. (We chose not to present the result initially in this generalized form to avoid confusing the reader with further notational difficulties which are easy to insert at a later point.)

Our proof of Theorem \ref{Theorem1} is largely self-contained and can replace the more complicated (but maybe conceptually easier)  proofs given in the first paper in the case $k=2$ and $k=3$.
We will first treat the case of pair correlation $k=2$ in detail in the next section. The following three sections will handle the case $k=3$ of triple correlation, where we give a more detailed proof then strictly necessary in order to provide examples for the general case. In the sixth section we prove the general case, and the following section discusses the computation of the constants $\mathcal{C}_k$.  In the final section we  obtain the generalized result mentioned above for use in some applications.

\emph{ Notation.} We retain the notation from our first paper which we will define here when it first appears. However, there are a few conventions we will use so frequently throughout the paper that they need to be mentioned immediately.  We will  use the Iverson notation which was defined below \eqref{1.8}. We will take $\epsilon$ to be any sufficiently small positive number whose value can be changed from equation to equation, and similarly $C$, $c$, and $c'$ will denote  small fixed constants whose value may change from equation to equation. For a vector $\text{\boldmath$j$}=(j_1,j_2,\cdots , j_k)$ we let $||\text{\boldmath$j$}|| = \max |j_i|$. A dash on a summation sign $\sum '$ indicates all the summation variables are relatively prime with each other, and further any sum without a lower bound on the summation variables will have the variables start with the value 1. Empty sums will have the value zero, and empty products will have the value 1.

\section*{ Acknowledgment}

The method we use to prove Theorem \ref{Theorem1} was suggested by Peter Sarnak during a talk given by the first-named author at the First Workshop on $L$-functions and Random Matrices at the American Institute of Mathematics in May 2001, and also suggested by John Friedlander after the talk.  Following the conference, the first-named author visited AIM where the method  was worked out jointly with Brian Conrey and David Farmer. Farmer has written a program in Mathematica to compute the constants $\mathcal{ C}_k$ and similar types of constants which has aided us in our work. The authors would like to thank these individuals and the American Institute of Mathematics.

\section{ Pair Correlation}

In this section we will prove the pair correlation case $k=2$ in Theorem \ref{Theorem1}. The proof is more difficult than our earlier proof in \cite{GYI}, but the method is self-contained and generalizes nicely to higher correlations. The result we prove in this section is stated in the following theorem.

\begin{theorem}\label{Theorem2.1} Let
\begin{equation}  \mathcal{S}_2(j) = \sum_{n=1}^N\Lambda_R(n)\Lambda_R(n+j).  \label{2.1} \end{equation}
Then with $\mathcal{D}$ a constant and $c'$ a small positive constant, we have
\begin{equation} \mathcal{S}_2(0) =  N\log R + \mathcal{D}N +  O( Ne^{-c'\sqrt{\log R}}) +O(R^2), \label{2.2}
\end{equation}
and, if $j\neq 0$, and letting $\text{\boldmath{$j$}}=(0,j)$, then for $\log |2j|\ll \log R$,
\begin{equation} \mathcal{S}_2(j) =  N\gs(\text{\boldmath{$j$}}) +  O( Ne^{-c'\sqrt{\log R}}) +O(R^2). \label{2.3}\end{equation}\end{theorem}

We have
\begin{equation}\begin{split} \mathcal{S}_2(j) &= \sum_{n=1}^N\Lambda_R(n)\Lambda_R(n+j) \\
&=\sum_{n\le N} \sum_{\substack{d,e\le R \\ d|n, \ e|n+j}} \mu(d)\mu(e) \log (R/d)\log (R/e)\\
&=\sum_{d,e \le R}\mu(d)\mu(e) \log (R/d)\log (R/e) \sum_{\substack{ n\le N\\ d|n \\ e|n+j}}1. \end{split} \label{2.4} \end{equation}
The divisibility conditions $d|n$ and $e|n+j$ can only be satisfied if $(d,e)|j$, in which case $n$ will run through a unique residue class modulo $[d,e]$. Hence
\begin{equation}
\sum_{\substack{ n\le N\\ d|n \\ e|n+j}}1 = [(d,e) | j]\left(\frac{N}{[d,e]} +O(1)\right),
\label{2.5} \end{equation}
and therefore
\begin{equation}
\mathcal{S}_2(j) = N \sum_{\substack{d,e \le R\\ (d,e)|j}}\frac{\mu(d)\mu(e)}{[d,e]} \log (R/d)\log (R/e) +O(R^2). \label{2.6}\end{equation}
To handle the least common multiple in the sum above, we let  $d = b_1 b_{12}$ and $e=b_2b_{12}$, where $b_{12} =(d,e)$. Then $[d,e]=b_1b_2b_{12}$ and $b_1$, $b_2$, and $b_{12}$ are pairwise relatively prime. (This notation conforms with the notation we will use in the general case.) Hence we have
\begin{equation}\begin{split} \mathcal{S}_2(j) & = N\sumprime_{\substack{b_1b_{12}\le R \\ b_2b_{12} \le R\\ b_{12} |j }} \frac{\mu(b_1)\mu(b_2)\mu^2(b_{12})}{b_1b_2b_{12}}\log\frac{R}{b_1b_{12}}\log\frac{R}{b_2b_{12}} +O(R^2) \\
& = NT_2(j) +O(R^2), \end{split}\label{2.7}\end{equation}
where the dash on the summation sign indicates the summation variables are pairwise relatively prime with each other.
At this point in the earlier proof in \cite{GYI} we summed over each variable individually. Here instead we apply the formula, for $m \ge 2$ and $c>0$,
\begin{equation} \frac{(m-1)!}{ 2\pi i}\int_{c-i\infty}^{c+i\infty} \frac{x^s}{ s^m}\, ds =
\left\{ \begin{array}{ll}
       0,
        &\text{if $0<x\le 1$,
    }  \\
(\log x)^{m-1},& \text{if $x\ge 1$.}
\end{array}
\right. \label{2.8}
\end{equation}
Letting $m=2$ in this formula, and denoting the vertical line contour $c+it$, $-\infty < t < \infty$  by $(c)$, we see that for $c_1 ,c_2 >0$,
\begin{equation}T_2(j)=
\frac {1}{ (2\pi i)^2}\mathop{\int }_{(c_2)\ }\! \mathop{\int}_{(c_1)\ } F(s_1,s_2)\frac{R^{s_1}}{ {s_1}^2}\frac{R^{s_2}}{ {s_2}^2}\, ds_1\,ds_2 ,\label{2.9}\end{equation}
where
\begin{equation}
F(s_1,s_2) = \sumprime_{\substack{1\le b_1,b_2,b_{12} < \infty\\ b_{12}|j}} \frac{\mu(b_1)\mu(b_2)\mu^2(b_{12})}{{b_1}^{1+s_1}{b_2}^{1+s_2}{b_{12}}^{1+s_1+s_2}}.\label{2.10}
\end{equation}
Letting $s_1 = \sigma_1 + it_1$ and $s_2=\sigma_2+it_2$, we see $F(s_1,s_2)$ is analytic in $s_1$ and $s_2$ for $ \sigma_1 >0$ and $\sigma_2 >0$. To analytically continue $F(s_1,s_2)$ to the left we start with the product representation
\begin{equation}
F(s_1,s_2) = \prod_{p \ndiv j}\left( 1 - \frac{1}{p^{1+s_1}} - \frac{1}{p^{1+s_2}}\right)\prod_{p | j}\left( 1 - \frac{1}{p^{1+s_1}} - \frac{1}{p^{1+s_2}}+ \frac{1}{p^{1+s_1+s_2}}\right). \label{2.11}
\end{equation}
There are two cases to consider. If $j=0$ then the first product is empty and $F(s_1,s_2)$ is equal to the second product over all primes, while if $j\neq 0$ then the second product is a finite product, and the analytic continuation of $F(s_1,s_2)$ depends on the first product. We will deal with each case separately.

\emph{Case 1.} Assume $j=0$. In this case, we have
\begin{equation}
F(s_1,s_2) = \prod_{p }\left( 1 - \frac{1}{p^{1+s_1}} - \frac{1}{p^{1+s_2}}+ \frac{1}{p^{1+s_1+s_2}}\right). \label{2.12}
\end{equation}

Since for $\mathrm{Re}(s)>1$
\begin{equation} \zeta(s) = \prod_p\left( 1 - \frac{1}{p^s}\right)^{-1}, \label{2.13}\end{equation}
and, for $\min(1+\sigma_1,1+\sigma_2,1+\sigma_1+\sigma_2) >0$,
\[\begin{split} 1 - \frac{1}{p^{1+s_1}} - \frac{1}{p^{1+s_2}} + \frac{1}{p^{1+s_1+s_2}} & = \left( 1 - \frac{1}{p^{1+s_1}}\right)\left( 1 - \frac{1}{p^{1+s_2}}\right) \left( 1 - \frac{1}{p^{1+s_1 + s_2}}\right)^{-1} \\ & \hskip .3in +O\left(\frac{1}{p^{2+\sigma_1+\sigma_2}}\Big(1 + \frac{1}{p^{\sigma_1}} + \frac{1}{p^{\sigma_2}}+ \frac{1}{p^{\sigma_1+ \sigma_2}}\Big)\right),\end{split} \]
we see that
\begin{equation} F(s_1,s_2) = \frac{\zeta(1+s_1+s_2)}{\zeta(1+s_1)\zeta(1+s_2)}h(s_1,s_2) \label{2.14}\end{equation}
where
\begin{equation} h(s_1,s_2) = \prod_p \frac{\left( 1 - \frac{1}{p^{1+s_1}} - \frac{1}{p^{1+s_2}} + \frac{1}{p^{1+s_1+s_2}}\right)}{ \left( 1 - \frac{1}{p^{1+s_1}}\right)\left( 1 - \frac{1}{p^{1+s_2}}\right)} \left( 1 - \frac{1}{p^{1+s_1 + s_2}}\right), \label{2.15}\end{equation}
and
\[ h(s_1,s_2) = \prod_{p}\left( 1 + O(\frac{1}{p^{2 + \nu}})\right) \]
with $\nu = (\sigma_1 + \sigma_2 )+ \min(0, \sigma_1,\sigma_2,\sigma_1+\sigma_2)$. Hence
we see $h(s_1,s_2)$ is analytic in $s_1$ and $s_2$ for $\sigma_1 > -\frac{1}{4}$ and $\sigma_2 > -\frac{1}{4}$. Therefore \eqref{2.14} provides the analytic continuation of $F(s_1,s_2)$ into this region and shows that  $|h(s_1,s_2)|\ll 1$ as $|t_1|, |t_2| \to \infty$. We also see that $h(0,0)=1$.

It should be kept in mind that while $F(s_1,s_2)$ is a function of two complex variables, we will at each stage in our calculation of \eqref{2.9} always fix one variable and treat $F$ as a function of the other single complex variable.  Thus, in the analysis that follows we assume in each equation that one of the complex variables $s_1$ or $s_2$ is fixed. When we treat the general $k$-correlations we will always fix all except one complex variable in $F$ at each stage in the calculation.

We now  isolate the dominant part of $F(s_1,s_2)$ around $s_1=0$, $s_2=0$. Since
\[(s-1)\zeta(s) = 1 + \sum_{n=1}^\infty a_ns^n\]
is analytic in the entire $s$-plane, we write
\begin{equation} \begin{split} F(s_1,s_2) &= \left(  \frac{(s_1+s_2)\zeta(1+s_1+s_2)}{s_1\zeta(1+s_1)s_2\zeta(1+s_2)}h(s_1,s_2)\right) \frac{s_1s_2}{s_1+s_2} \\& = k(s_1,s_2) \frac{s_1s_2}{s_1+s_2}. \end{split}\label{2.16}\end{equation}
We see
\begin{equation}  k(s_1,0)=k(0,s_2)=1 , \label{2.17} \end{equation}
and $k(s_1,s_2)$ is analytic everywhere in $s_1$ and $s_2$ except for
poles at the zeros of $\zeta(1+s_1)$ and $\zeta(1+s_2)$. To avoid these poles, we use a classical zero-free region result.  By  Theorem 3.11 and (3.14.8) of \cite{T}  there exists a small positive constant $C$ such that $\zeta(\sigma +it)\neq 0$ in the region
\begin{equation} \sigma \ge 1 -\frac{C}{ \log(|t|+2)} \label{2.18}\end{equation}
 for all $t$, and further
\begin{equation} \zeta(\sigma +it)-\frac{1}{\sigma -1 +it} \ll \log(|t|+2), \hskip .4in
\frac{1}{ \zeta(\sigma +it)} \ll \log(|t|+2)
\label{2.19}
\end{equation}
in this region. (There are stronger results but this suffices for our needs.)
Hence $k(s_1,s_2)$ is analytic in $s_1$ and $s_2$ in the region
\begin{equation} \sigma_1 \ge -\frac{C}{\log (|t_1|+2)}, \hskip .4in \sigma_2 \ge -\frac{C}{\log (|t_2|+2)}, \label{2.20}\end{equation}
and we also see by \eqref{2.19} that in this region
\begin{equation}  F(s_1,s_2)\ll \log(|t_1|+2)\log(|t_2|+2)\max\big( \frac{1}{|s_1+s_2|}, \log( |t_1+t_2|+2)\big)  .\label{2.21} \end{equation}

Returning to \eqref{2.9} in this case,  we have
\begin{equation}  T_2(0) = \frac {1}{ (2\pi i)^2}\mathop{\int }_{(c_2)\ }\! \mathop{\int}_{(c_1)\ } k(s_1,s_2) \frac{R^{s_1+s_2}}{s_1s_2(s_1+s_2)} \, ds_1\,ds_2 , \label{2.22}
\end{equation}
where we now take $c_2 >c_1 >0$. To evaluate this integral, we move the contours to the edge of the region \eqref{2.20} and evaluate the residues. Let $\mathcal{ L}$ denote the contour given by
\begin{equation} s= -\frac{C}{ \log(|t|+2)} +it.\label{2.23}\end{equation}
and let $\mathcal{L}_j$ denote the same contour except with $C=C_j$. Since we want to avoid integrating over the pole $s_1=-s_2$ of the integrand on the contour we move to,  we will move $(c_j)$ to $\mathcal{L}_j$ and choose $C_2 = \frac{1}{2}C_1$ so that the contours $\mathcal{L}_1$ and $\mathcal{L}_2$ are well separated with
\[|s_1+s_2|\gg \frac{1}{\log(|s_1| +2)}\]
 on these contours. Now,
provided $|s_1|\gg 1$, $|s_2|\gg 1$, and $|s_1+s_2|\gg \frac{1}{\log(|s_1|+2)}$,
\begin{equation} \left|F(s_1,s_2)\frac{R^{s_1+s_2}}{{s_1}^2{s_2}^2}\right| \ll R^{\sigma_1+\sigma_2}\frac{ \log^2(|t_1|+2)\log^2(|t_2|+2)}{ (|t_1|+2)^2(|t_2|+2)^2} .\label{2.24}\end{equation}
 in the region \eqref{2.20}.
This bound shows the integral in \eqref{2.22} is absolutely convergent, and, since we can easily arrange $(c_1)$ and $(c_2)$ to meet the conditions necessary for this bound, which are also satisfied on $\mathcal{L}_1$ and $\mathcal{L}_2$, justifies the residue calculations that follow. Moving first $(c_1)$ to $\mathcal{L}_1$, we encounter a simple pole at $s_1=0$ with residue $ \frac{R^{s_2}}{{s_2}^2}$
by \eqref{2.17}, and hence we obtain using \eqref{2.8}
\[ \begin{split} T_2(0) & = \frac {1}{ 2\pi i}\mathop{\int }_{(c_2)\ }  \frac{R^{s_2}}{{s_2}^2}\, ds_2+
\frac {1}{ (2\pi i)^2}\mathop{\int }_{(c_2)\ }\! \mathop{\int}_{\mathcal{L}_1\ } k(s_1,s_2) \frac{R^{s_1+s_2}}{s_1s_2(s_1+s_2)} \, ds_1\,ds_2 \\ & = \log R + \frac {1}{ (2\pi i)^2}\mathop{\int }_{(c_2)\ }\! \mathop{\int}_{\mathcal{L}_1\ } k(s_1,s_2) \frac{R^{s_1+s_2}}{s_1s_2(s_1+s_2)} \, ds_1\,ds_2 . \end{split}\]
In the double integral above we move the contour $(c_2)$ to the left to $\mathcal{L}_2$ and encounter two simple poles at $s_2=0$ and $s_2=-s_1$,   with residues
\[ \frac{R^{s_1}}{{s_1}^2}, \quad \mathrm{and} \  - \frac{k(s_1,-s_1)}{{s_1}^2}, \]
respectively.
Hence
\begin{equation}  T_2(0) = \log R +\mathcal{D} + \frac {1}{ 2\pi i}\mathop{\int }_{\mathcal{L}_1\ } \frac{R^{s_1}}{{s_1}^2}\, ds_1  +
\frac {1}{ (2\pi i)^2}\mathop{\int }_{\mathcal{L}_2\ }\! \mathop{\int}_{\mathcal{L}_1\ } F(s_1,s_2) \frac{R^{s_1+s_2}}{{s_1}^2{s_2}^2} \, ds_1\,ds_2 , \label{2.25}\end{equation}
where
\begin{equation} \mathcal{D} = - \frac {1}{ 2\pi i}\mathop{\int }_{\mathcal{L}_1\ }   \frac{k(s_1,-s_1)}{{s_1}^2}\, ds_1 . \label{2.26} \end{equation}
The first integral in \eqref{2.25} is zero since we may move the contour to the left without encountering any singularities to $-\infty$ where it vanishes.  Using \eqref{2.24} the remaining integral in \eqref{2.25} is bounded by
\[  \ll \left( \int_{-\infty}^\infty R^{-\frac{c}{ \log(|t|+2)}}\frac{\log^2(|t|+2)}{{t}^2+1}\, dt\right)^2 .\]
The integral here is, for any $w\ge 2$,
\[ \ll \log^2 w \int_0^w R^{-\frac{c}{\log(|t|+2)}}\, dt + \int_w^\infty \frac{\log^2 t}{ t^2}\, dt
 \ll w(\log^2 w) e^{\frac{-c\log R}{ \log w}} + \frac{\log^2 w }{ w} ,\]
and on choosing $\log w = \frac{1}{ 2}\sqrt{c\log R}$ this is
\[ \ll  (c\log R) e^{-\frac{1}{ 2} \sqrt{c\log R}} \ll e^{-c'\sqrt{\log R}}.\]
We note for future use that the above argument shows that with $s=\sigma + it$, $\sigma = -\frac{c}{\log(|t|+2)}$ and $B$ a constant
\begin{equation} \int_{-\infty}^{\infty}R^\sigma \frac {\log^B(|s|+2)}{|s|^2} \, dt \ll e^{-c'\sqrt{\log R}}, \label{2.27} \end{equation}
since here
\[ \frac{1}{|s|^2}\ll \frac{\log(|s|+2)}{2+|t|^2}.\]

We conclude that
\begin{equation}  T_2(0) = \log R + \mathcal{D} + O(e^{-c'\sqrt{\log R}}),\label{2.28}
\end{equation}
which by \eqref{2.7} proves the first part of Theorem \ref{Theorem2.1}.

For computation purposes we note that the main term in the above analysis is obtained directly from
\begin{equation} \frac {1}{ (2\pi i)^2}\mathop{\int }_{(c_2)\ }\! \mathop{\int}_{(c_1)\ } \frac{R^{s_1+s_2}}{s_1 s_2 (s_1+s_2)} \, ds_1\,ds_2 , \label{2.29}\end{equation}
since the main term in \eqref{2.22} occurs from the residues with $k(0,0)$ as a factor.
We  make  the change of variables $w_1 = (\log R) s_1$ and  $w_2 = (\log R)s_2$. Since we can take $c_j/\log R$ in place of $c_j$  for $j=1,2$ in the double integral, we conclude that the main term is
\begin{equation}   \mathcal{C}_2(2) \log R ,\label{2.30}
\end{equation}
where, on returning to the variables $s_1$ and $s_2$,
\begin{equation} \mathcal{C}_2(2) = \frac {1}{ (2\pi i)^2}\mathop{\int }_{(c_2)\ }\! \mathop{\int}_{(c_1)\ } \frac{e^{s_1+s_2}}{s_1s_2(s_1+s_2)} \, ds_1\,ds_2 .\label{2.31} \end{equation}
We will discuss this integral and the computation of $\mathcal{C}_k(k)$ in Section 7.

\emph{Case 2.} $j\neq 0$. In this case, the first product in \eqref{2.11} is the dominant one, and therefore
\begin{equation} F(s_1,s_2) =  \frac{1}{\zeta(1+s_1)\zeta(1+s_2)}h_j(s_1,s_2) \label{2.32} \end{equation}
where
\begin{equation}
h_j(s_1,s_2) = \prod_{p \ndiv j}\left( \frac{1 - \frac{1}{p^{1+s_1}} - \frac{1}{p^{1+s_2}}}{\big( 1 - \frac{1}{p^{1+s_1}}\big)\big( 1 - \frac{1}{p^{1+s_2}}\big)}\right)\prod_{p | j}\left( \frac{1 - \frac{1}{p^{1+s_1}} - \frac{1}{p^{1+s_2}}+ \frac{1}{p^{1+s_1+s_2}}}{\big( 1 - \frac{1}{p^{1+s_1}}\big)\big( 1 - \frac{1}{p^{1+s_2}}\big)}\right),\label{2.33}
\end{equation}
and we see $h_j(s_1,s_2)$ is analytic in the region $\min(1+ \sigma_1 +\sigma_2, 1+\sigma_1 , 1+\sigma_2) > 0$, or in particular if $\sigma_1 ,\sigma_2 > -\frac{1}{2}$.  Hence $F(s_1,s_2)$ is analytic in the region \eqref{2.20}. We also see
\[  h_j(0,0) = \prod_{ p \ndiv j}\Big( \frac{1-\frac{2}{p}}{\big(1- \frac{1}{p}\big)^2} \Big) \prod_{ p | j} \Big( \frac{ 1}{1-\frac{1}{p}}\Big)  = \prod_p\Big(1 - \frac{1}{p}\Big)^{-2}\Big( 1 - \frac{\nu_p(j)}{p}\Big),  \]
where, since $j\neq 0$, $\nu_p(j)$ is 1 if $p|j$ and $\nu_p(j) =2$ if $p\ndiv j$. Thus $\nu_p(j)= \nu_p(\text{\boldmath{$j$}})$ when $\text{\boldmath{$j$}} =(0,j)$, and we have shown that
\begin{equation} h_j(0,0) = \gs(\text{\boldmath{$j$}}).\label{2.34} \end{equation}
We now move the contour integrals in \eqref{2.9} to the left as before and encounter simple poles at $s_1=0$ and $s_2=0$ which give
\begin{equation} T_2(j) = \gs(\text{\boldmath{$j$}})  + \frac {1}{ (2\pi i)^2}\mathop{\int}_{\mathcal{L}_2\ } \! \mathop{\int}_{\mathcal{L}_1\ }F(s_1,s_2)\frac{R^{s_1+s_2}}{{s_1}^2{s_2}^2} \, ds_1\,
ds_2 . \label{2.35} \end{equation}
By the estimate preceding (3.5) of \cite{GYI} we have for $\delta >0$
\begin{equation} \prod_{p|j}\left(1 + O(\frac{1}{p^{1-\delta}})\right) \ll \exp( D(\log |2j|)^{\delta}) \label{2.36}\end{equation}
for some positive constant $D$,
and therefore in the  $-\frac{1}{2} < \sigma_1 ,\sigma_2 \le 0$ we see that
\[h_j(s_1,s_2) \ll \prod_{p|j}\Big( 1 + O\big(\frac{1}{p^{1+\sigma_1+\sigma_2}}\big)\Big) \ll \exp(D(\log |2j|)^{-\sigma_1 -\sigma_2}). \]
We conclude by \eqref{2.19} and \eqref{2.27} that the double integral in \eqref{2.35} is
\[ \ll \left( \exp( D(\log |2j|)^{\frac{C_1+C_2}{\log 2}})\int_{-\infty}^\infty \frac{ R^{\sigma}\log^2( |t|+2)}{|s|^2} \, dt\right)^2 \ll e^{-c'\sqrt{\log R}},\]
on taking (as we may) $C_1 +C_2 < \frac{1}{2}\log 2$ and using $\log|2j|\ll \log R$.  This completes the proof of Theorem \ref{Theorem2.1}.

\section{Initial treatment of triple correlation}
We now consider the triple correlations of $\Lambda_R(n)$. The procedure is the same as for the pair correlations only the calculations are more complicated. We have been more detailed then necessary in order to provide motivation for the general case. We will prove the following theorem.
\begin{theorem} \label{Theorem3.1} For constants $\mathcal{D}_0$ and $\mathcal{D}_1$ we have
\begin{equation} S_3(0) = \sum_{n=1}^N \Lambda_R(n)^3 = \frac{3}{4}N\log^2 R + \mathcal{D}_1 N\log R + \mathcal{D}_0 N +O(Ne^{-c'\sqrt{\log R}}) + O(R^3). \label{3.1}\end{equation}
For $\text{\boldmath$j$}= (j_1,j_2)$, $j_1\neq j_2 \neq 0$, with $\log 2||\boldsymbol{j}|| \ll \log R$, and $\text{\boldmath$a$}=(2,1)$ we have
\begin{equation} \mathcal{S}_3(N,\text{\boldmath$j$},\text{\boldmath$a$}) = \gs(\text{\boldmath$j$})N \log R + \mathcal{D}_2(\text{\boldmath$j$})N + O(N e^{-c'\sqrt{\log R}}) + O(R^3) ,\label{3.2}\end{equation}
where $\mathcal{D}_2(\text{\boldmath$j$})$ is an arithmetic function which satisfies the bound $\mathcal{D}_2(\text{\boldmath$j$}) \ll \exp(||\boldsymbol{j}||^\delta)$ for any fixed $\delta>0$. For $\text{\boldmath$j$}= (j_1,j_2,j_3)$, $j_1\neq j_2 \neq j_3 \neq 0$ with $\log 2||\boldsymbol{j}|| \ll \log R$, and $\text{\boldmath$a$}=(1,1,1)$, we have
\begin{equation} \mathcal{S}_3(N,\text{\boldmath$j$},\text{\boldmath$a$}) = \gs(\text{\boldmath$j$})N + O( N e^{-c'\sqrt{\log R}}) + O(R^3) .\label{3.3}\end{equation}
\end{theorem}

Let
\begin{equation}\mathcal{ S}_3(j_1,j_2,j_3) = \sum_{n=1}^N \Lambda_R(n+j_1)\Lambda_R(n+j_2)\Lambda_R(n+j_3).
\label{3.4}
\end{equation}
Expanding,  we have
\[\mathcal{ S}_3(j_1,j_2,j_3) = \sum_{d_1,d_2,d_3\le R}
\mu(d_1)\log (R/d_1)\mu(d_2)\log (R/d_2)
\mu(d_3)\log (R/d_3)\sum_{\stackrel{\stackrel{\stackrel{\scr n\le N}
 {\scr d_1|n+j_1 }}  {\scr d_2|n+j_2 }} { \scr d_3|n+j_3}}1.\]
The sum over $n$ is zero unless $(d_1,d_2)|j_2-j_1$, $(d_1,d_3)|j_3-j_1$, and $(d_2,d_3)|j_3-j_2$, in which case the sum runs through a unique residue class modulo $[d_1,d_2,d_3]$, and we have
\[\sum_{\stackrel{\stackrel{\stackrel{\scr n\le N}
 {\scr d_1|n+j_1 }}  {\scr d_2|n+j_2 }} { \scr d_3|n+j_3}}1 = \frac{N}{[d_1,d_2,d_3]} + O(1).\]
We conclude, letting
\begin{equation} L_i(R) = \log \frac{R}{d_i}, \label{3.5}\end{equation}
that
\begin{equation}
\begin{array}{lcl}\mathcal{ S}_3(j_1,j_2,j_3) &= & {\displaystyle N\sum_{\stackrel{\stackrel{\stackrel{\scr d_1,d_2,d_3\le R} {\scr (d_1,d_2)|j_2-j_1 }}  {\scr  (d_1,d_3)|j_3-j_1} }{ \scr (d_2,d_3)|j_3-j_2}}
\frac{\mu(d_1)\mu(d_2)\mu(d_3)}{ [d_1,d_2,d_3] }L_1(R)L_2(R)L_3(R) + O(R^3) }\\
&=& N T_3(j_1,j_2,j_3) +O(R^3).
\end{array}
 \label{3.6}
\end{equation}
We now decompose $d_1$, $d_2$, and $d_3$ into relatively prime factors
\begin{eqnarray*} d_1&=&b_1b_{12}b_{13}b_{123} ,\\
d_2&=&b_2b_{12}b_{23}b_{123}, \\
d_3&=&b_3b_{13}b_{23}b_{123} ,
\end{eqnarray*}
where $b_\chi$ is a divisor of the $d_i$'s where $i$ occurs in the digits of $\chi$. Since the $d_i$'s are squarefree,  these new variables are pairwise relatively prime.  We will let $\mathcal{ D}$ denote the set of $b_\chi$'s which satisfy the conditions
\[ b_{12}b_{123}|j_2-j_1, \quad b_{13}b_{123}|j_3-j_1, \quad b_{23}b_{123}|j_3-j_2 .\]
Then we have
\begin{equation} \begin{split}
 & T_3(j_1,j_2,j_3) = \\ &\sumprime_{\substack{b_1b_{12}b_{13}b_{123}\le R\\ b_2b_{12}b_{23}b_{123}\le R  \\
b_3b_{13}b_{23}b_{123}\le R \\  \mathcal{D}}}\frac{ \mu(b_1)\mu(b_2)\mu(b_3)\mu^2(b_{12})\mu^2(b_{13})\mu^2(b_{23})\mu(b_{123}) }{ b_1 b_2 b_3 b_{12}b_{13}b_{23}b_{123}} L_1(R)L_2(R)L_3(R)
 .\label{3.7}\end{split}
\end{equation}
 Now by \eqref{2.8} we have
\begin{equation}T_3(j_1,j_2,j_3)=
\frac {1}{ (2\pi i)^3} \mathop{\int}_{(c_3)\ }\!\mathop{\int}_{(c_2)\ } \!\mathop{\int}_{(c_1)\ } F(s_1,s_2,s_3)\frac{R^{s_1}}{ {s_1}^2}\frac{R^{s_2}}{ {s_2}^2}\frac{R^{s_3}}{ {s_3}^2}\, ds_1\,ds_2 \, ds_3 ,\label{3.8}\end{equation}
where
\begin{equation}
F(s_1,s_2,s_3) = \sumprime_{\mathcal{D}} \frac{\mu(b_1)\mu(b_2)\mu(b_3)\mu^2(b_{12})\mu^2(b_{13})\mu^2(b_{23})\mu(b_{123})}{{b_1}^{1+s_1}{b_2}^{1+s_2}{b_3}^{1+s_3}{b_{12}}^{1+s_1+s_2}{b_{13}}^{1+s_1+s_3}{b_{23}}^{1+s_2+s_3}{b_{123}}^{1+s_1+s_2+s_3}}.\label{3.9}
\end{equation}
We define $J_{12}$, $J_{13}$, and $J_{23}$ to be the largest squarefree divisor of $j_2 -j_1$, $j_3-j_1$, and $j_3-j_2$ respectively. Let
\begin{equation} \Delta_3 = (j_2 -j_1)(j_3-j_1)(j_3-j_2), \qquad \kappa = (J_{12},J_{13},J_{23}).\label{3.10} \end{equation}
Since
$j_3-j_2 = (j_3-j_1) -(j_2-j_1)$, we see that if a prime $p \ndiv \kappa$, then it can divide at most one of numbers $J_{12}$, $J_{13}$, $J_{23}$.
Here for convenience we will extend the usual definition of the gcd by $(a,0) = a$, and $(0,0) =0$.
With $s_n = \sigma_n +it_n$, $1\le n\le 3$, we have for $\sigma_n >0$ with $1\le n \le 3$ the product representation
\begin{equation}
F(s_1,s_2,s_3) = h_1(\kappa)h_2(\kappa, J_{12})h_2(\kappa, J_{13})h_2(\kappa, J_{23})h_3(\Delta_3) , \label{3.11}
\end{equation}
where
\begin{equation}\begin{split}
&h_1(\kappa) = h_1(\kappa, s_1,s_2,s_3) =\\& \prod_{p | \kappa}\left( 1 - \frac{1}{p^{1+s_1}} - \frac{1}{p^{1+s_2}}-\frac{1}{p^{1+s_3} }+\frac{1}{p^{1+s_1+s_2}}+\frac{1}{p^{1+s_1+s_3}}+ \frac{1}{p^{1+s_2+ s_3}}-\frac{1}{p^{1+s_1+s_2+s_3}}\right),\end{split}\label{3.12}
\end{equation}
\begin{equation}
h_{2}(\kappa , J_{ij}) = h_{2}(\kappa , J_{ij},s_1,s_2,s_3) = \prod_{\substack{p\ndiv \kappa\\ p | J_{ij}}}\left( 1 - \frac{1}{p^{1+s_1}} - \frac{1}{p^{1+s_2}}-\frac{1}{p^{1+s_3} }+\frac{1}{p^{1+s_i+s_j}}\right),\label{3.13}\end{equation}
for $1\le i<j\le 3$, and
\begin{equation}
h_3(\Delta_3) =  h_3(\Delta_3, s_1,s_2,s_3)= \prod_{p \ndiv \Delta_3}\left( 1 - \frac{1}{p^{1+s_1}} - \frac{1}{p^{1+s_2}}-\frac{1}{p^{1+s_3} }\right).\label{3.14}\end{equation}

There are three cases to consider, when $j_1 = j_2 =j_3$, when  $j_1 = j_2 \neq j_3$ (or a permutation of this), and when $j_1 \neq j_2 \neq j_3$. These cases correspond to the three numbers $j_2-j_1$, $j_3-j_1$, and $j_3-j_2$ falling into 1, 2, or 3 different residue classes for all sufficiently large primes. We consider each case separately.

\section{ $S_3(0)$}
We consider the first case when $j_1=j_2=j_3$, where \eqref{3.11} becomes
\begin{equation} \begin{split}
&F(s_1,s_2,s_3)= h_1(0,s_1,s_2,s_3)=\\& \prod_{p }\left( 1 - \frac{1}{p^{1+s_1}} - \frac{1}{p^{1+s_2}}-\frac{1}{p^{1+s_3} }+\frac{1}{p^{1+s_1+s_2}}+\frac{1}{p^{1+s_1+s_3}}+ \frac{1}{p^{1+s_2+ s_3}}-\frac{1}{p^{1+s_1+s_2+s_3}}\right),\end{split}\label{4.1}
\end{equation}
and we see
\begin{equation} F(s_1,s_2,s_3) = \frac{\zeta(1+s_1+s_2 )\zeta(1+s_1 + s_3)\zeta(1+ s_2 + s_3) }{\zeta(1+s_1) \zeta(1+s_2)\zeta(1+s_3)\zeta(1+s_1+s_2 + s_3)}h(s_1,s_2,s_3), \label{4.2}\end{equation}
where $h(0,0,0)=1$, and
\begin{equation} h(s_1,s_2,s_3) = \prod_p\left( 1 +O(\frac{1}{p^{2+\nu}})\right), \label{4.3}\end{equation}
for  $\sigma_1 +\sigma_2 +\sigma_3 > -\frac{1}{2} $, where, assuming $\sigma_1\le \sigma_2 \le \sigma_3$,
\[ \nu \le  \sigma_1 +\sigma_2 + \min( 0, \sigma_1, \sigma_1+\sigma_2, \sigma_1 +\sigma_2 + \sigma_3).\]
Hence $h(s_1,s_2,s_3)$ is analytic in $s_1$, $s_2$, and $s_3$ for $\sigma_1,\sigma_2, \sigma_3 > -\frac{1}{5}$, and $|h(s_1,s_2,s_3)|\ll 1$ as $|t_1|,|t_2|,|t_3| \to \infty$.
Next, as in \eqref{2.16} we define $k(s_1,s_2,s_3)$ by
\begin{equation}
F(s_1,s_2,s_3) = k(s_1,s_2,s_3)\frac{s_1s_2s_3(s_1+s_2+s_3)}{(s_1+s_2)(s_1+s_3)(s_2+s_3)}.  \label{4.4}\end{equation}
Hence
\begin{equation} k(s_1,s_2,0) =k(s_1,0,s_3)=k(0,s_2,s_3)=1 \label{4.5} \end{equation}
and as before $k(s_1,s_2,s_3)$  is analytic in the region
\begin{equation} \sigma_j \ge -\frac{C}{\log(|t_j|+2)}, \quad 1\le j\le 3, \label{4.6} \end{equation}
for a small enough constant $C$, and in this region by \eqref{2.19}
\begin{equation}  F(s_1,s_2,s_3) \ll \left(\prod_{i=1}^3\log^2(|t_i|+2)\right)\left( \prod_{1\le i<j\le 3}\max\big(\frac{1}{|s_i+s_j|}, \log(|t_i| + 2)\big)\right). \label{4.7} \end{equation}

As before we let $\mathcal{ L}_j$ denote the contour from \eqref{2.23} with the constant $C=C_j$. Starting with a constant $C_1$ for which \eqref{2.18} holds, we take the constants $C_{j} = \frac{1}{2}C_{j-1}$  for the contour $\mathcal{L}_{j}$, $j=2,3$,  which arranges the contours so that $\mathcal{L}_j$ is to the left of $\mathcal{L}_{j+1}$ and these contours are well-spaced from each other. We let
\begin{equation} \mathcal{J}(c_1,c_2,c_3) = \frac {1}{ (2\pi i)^3}\mathop{\int}_{(c_3)\ }\!\mathop{\int}_{(c_2)\ } \!\mathop{\int}_{(c_1)\ }
k(s_1,s_2,s_3)\frac{R^{s_1+s_2+s_3}(s_1+s_2+s_3)}{s_1s_2s_3(s_1+s_2)(s_1+s_3)(s_2+s_3)}\, ds_1ds_2ds_3 .\label{4.8} \end{equation}
and see by \eqref{3.8} that for $c_3 > c_2>c_1 >0$ (or we can just take $c_j=j$)
\begin{equation} T_3(0) = \mathcal{J}(c_1,c_2,c_3).\label{4.9} \end{equation}
We first move $c_1$ to the left to $\mathcal{L}_1$, and encounter only the simple pole at $s_1=0$ since the other two poles of the integrand at $s_1=-s_2$ and $s_1=-s_3$ are to the left of $\mathcal{L}_1$. By \eqref{4.5} the residue of this pole at $s_1=0$ is
\[   \frac{R^{s_2+s_3}}{{s_2}^2{s_3}^2}, \] and therefore by \eqref{2.8}
\[ T_3(0) = \log^2 R + \mathcal{J}(\mathcal{L}_1,c_2,c_3). \]
Moving next $(c_2)$ to the left to $\mathcal{L}_2$, we encounter simple poles as $s_2=0$ and  $s_2=-s_1$ (but not $s_2=-s_3$) with residues
\[ \frac{R^{s_1+ s_3}}{{s_1}^2{s_3}^2}, \quad \mathrm{ and} \  - \frac{k(s_1,-s_1,s_3)R^{s_3}}{{s_1}^2(s_1+s_3)(-s_1+s_3)},\]
respectively,
and hence
\begin{equation} \begin{split}
T_3(0) &= \log^2 R +  \frac {1}{ (2\pi i)^2}\mathop{\int}_{\mathcal{L}_1\ }\!\mathop{\int}_{(c_3)\ }  \frac{R^{s_1+s_3}}{{s_1}^2{s_3}^2} \, ds_3\,ds_1   \\ & \qquad - \frac {1}{ (2\pi i)^2}\mathop{\int}_{\mathcal{L}_1\ }\! \mathop{\int}_{(c_3)\ }\!\frac{k(s_1,-s_1,s_3)R^{s_3}}{{s_1}^2 (s_1+s_3)(-s_1+s_3)} \, ds_3\,
 ds_1  +\mathcal{J}(\mathcal{L}_1,\mathcal{L}_2,c_3) .\end{split} \label{4.10}\end{equation}
The first integral above is zero since we can move $\mathcal{L}_1$ to the left to $-\infty$ where it vanishes without encountering any singularities. In the remaining terms we move $(c_3)$ to $\mathcal{L}_3$; the second integral encounters a simple pole at $s_3=-s_1$ but not the pole at $s_3 = s_1$ due to how we arranged the contours $\mathcal{L}_i$, and in $\mathcal{J}$ we pass simple poles $s_3=0$, $s_3=-s_1$, and  $s_3=-s_2$.  We conclude
\begin{equation} \begin{split}
T_3(0) &= \log^2 R  +  \frac {1}{ 2\pi i}\mathop{\int}_{\mathcal{L}_1\ }  \frac{k(s_1,-s_1,-s_1)R^{-s_1}}{2{s_1}^3} \,ds_1 \\& \qquad - \frac {1}{ (2\pi i)^2}\mathop{\int}_{\mathcal{L}_1\ }\!\mathop{\int}_{\mathcal{L}_3\ } \frac{k(s_1,-s_1,s_3)R^{s_3}}{{s_1}^2 (s_1+s_3)(-s_1+s_3)}   \, ds_3\,ds_1     + \frac {1}{ (2\pi i)^2}\mathop{\int}_{\mathcal{L}_1\ }\! \mathop{\int}_{\mathcal{L}_2\ }\!\frac{R^{s_1+s_2}}{{s_1}^2 {s_2}^2} \, ds_2\,
 ds_1 \\ & \qquad - \frac {1}{ (2\pi i)^2}\mathop{\int}_{\mathcal{L}_1\ }\!\mathop{\int}_{\mathcal{L}_2\ } \left(\frac{k(s_1,s_2,-s_1)R^{s_2}}{{s_1}^2 (s_1+s_2)(-s_1+s_2)}   +\frac{k(s_1,s_2,-s_2)R^{s_1}}{{s_2}^2 (s_1+s_2)(s_1-s_2)}   \right)\, ds_2\,ds_1 \\ & \qquad +\mathcal{J}(\mathcal{L}_1,\mathcal{L}_2,\mathcal{L}_3) \\&
= \log^2 R +  \frac {1}{ 2\pi i}\mathop{\int}_{\mathcal{L}_1\ }  \frac{k(s_1,-s_1,-s_1)R^{-s_1}}{2{s_1}^3} \,ds_1 -I_1 +I_2 - I_3 + \mathcal{J}(\mathcal{L}_1,\mathcal{L}_2,\mathcal{L}_3) .\end{split}\label{4.11} \end{equation}
The integral $I_2$ is zero since we may move the contours to the left to $-\infty$ where the integral vanishes. To estimate the remaining integrals, we note for $s_i$ on $\mathcal{L}_i$ we have  $|s_i|\gg \frac{1}{\log(|s_i| +2)}$, and
\[  |s_i \pm s_j| \gg \frac{1}{\log(|s_i| +2)}, \quad i\neq j. \]
Thus by \eqref{4.7} and \eqref{2.27}
\begin{equation} \begin{split} \mathcal{J}(\mathcal{L}_1,\mathcal{L}_2,\mathcal{L}_3) & \ll
\left( \int_{\mathcal{L} }R^{\sigma} \frac {\log^B(|s|+2)}{|s|^2}\, dt\right)^3 \\ &
\ll e^{-c'\sqrt{\log R}} . \end{split}\label{4.12} \end{equation}
Similarly we have in $I_1$
\[ \begin{split} \frac{k(s_1,-s_1,s_3)R^{s_3}}{{s_1}^2 (s_1+s_3)(-s_1+s_3)} & = - \frac{\zeta(1+s_1 +s_3)\zeta(1-s_1+s_3)R^{s_3}}{{s_1}^4{s_3}^2\zeta(1+s_1)\zeta(1-s_1)\zeta^2(1+s_3)}h(s_1,-s_1,-s_3)\\&
\ll \frac{R^{\sigma_3}}{|s_1|^4|s_3|^2}\log^4(|t_1|+2)\log^4(|t_3|+2)\\&
\ll\frac{R^{\sigma_3}}{|s_1|^2|s_3|^2}\log^6(|t_1|+2)\log^4(|t_3|+2)\end{split}\]
and thus \eqref{2.27} applies, and the same argument holds for $I_3$. Finally, for the first integral in \eqref{4.11}
we move the contour back to the right across the triple pole at $s_1=0$ to the contour ${\mathcal{L}_1}^*$ given by $s_1 = \frac{C_1}{\log(|t_1|+2)} +it$ which is the reflection of $\mathcal{L}_1$ over the imaginary axis. Thus this integral is
\begin{equation}  = \mathop{\mathrm{res}}_{s_1=0}\left(  \frac{k(s_1,-s_1,-s_1)R^{-s_1}}{2{s_1}^3}\right) + \frac{1}{2 \pi i}\mathop{\int}_{{\mathcal{L}_1}^* \ } \frac{k(s_1,-s_1,-s_1)R^{-s_1}}{2{s_1}^3}\, ds_1 . \label{4.13}\end{equation}
The integral here can now be estimated as before by \eqref{2.27} (because $1/|s_1|^3\ll \log(|t_1|+2)/|s_1|^2$ on the contour) with the same error term.
Finally, the triple pole has residue
\begin{equation}  \mathop{\mathrm{res}}_{s_1=0}\left( R^{-s_1}\frac{\zeta(1-2s_1)h(s_1,-s_1,-s_1}{{s_1}^6\zeta(1+s_1)\zeta^3(1-s_1)}\right) = -\frac{1}{4}\log^2R + \mathcal{D}_1 \log R +\mathcal{D}_0, \label{4.14}\end{equation}
for constants $\mathcal{D}_0$ and $\mathcal{D}_1$.
Combining our results, we conclude
\begin{equation} T_3(0) = \frac{3}{4}\log^2 R + \mathcal{D}_1 \log R +\mathcal{D}_0 +O(e^{-c'\sqrt{\log R}}). \label{4.15} \end{equation}
 This completes the proof  of \eqref{3.1} in Theorem \ref{3.1}.

For computational purposes we note that the main term is obtained from terms with $k(0,0,0)$ which from \eqref{4.8} is
\[ \frac {1}{ (2\pi i)^3}\mathop{\int}_{(c_3)\ }\!\mathop{\int}_{(c_2)\ } \!\mathop{\int}_{(c_1)\ } \frac{(s_1+s_2+s_3)R^{s_1+s_2+s_3}}{s_1s_2s_3(s_1+s_2)(s_1+s_3)(s_2+s_3)} \, ds_1\,ds_2\,ds_3.  \]
On making the change of variables $w_1 = (\log R) s_1$,  $w_2 = (\log R)s_2$,  $w_3 = (\log R)s_3$ and taking $c_j/\log R$ in place of $c_j$  for $1\le j\le 3$ this term becomes
\[ \mathcal{C}_3(3)\log^2R ,\]
where, on replacing $w_i$ by $s_i$,
\begin{equation} \mathcal{C}_3(3) = \frac {1}{ (2\pi i)^3}\mathop{\int}_{(c_3)\ }\!\mathop{\int}_{(c_2)\ } \!\mathop{\int}_{(c_1)\ } \frac{(s_1+s_2+s_3)e^{s_1+s_2+s_3}}{s_1s_2s_3(s_1+s_2)(s_1+s_3)(s_2+s_3)} \, ds_1\,ds_2\, ds_3 .\label{4.16} \end{equation}

\section{$\mathcal{ S}_3(N,\text{\boldmath$j$}, \text{\boldmath $a$})$ for $\text{\boldmath $a$}=(2,1)$ and $\text{\boldmath$a$}=(1,1,1)$}
We first consider the case $j_1=j_2 \neq j_3$. Let $j = j_3-j_1 =j_3-j_2$, and take $J$ to be the largest squarefree divisor of $j$. Then in \eqref{3.10} we have $\Delta_3 =0$ and $\kappa = J$, and \eqref{3.11} becomes
\begin{equation} \begin{split} F(s_1,s_2,s_3)&= h_1(J,s_1,s_2,s_3)\prod_{p\ndiv J}\left( 1 - \frac{1}{p^{1+s_1}} - \frac{1}{p^{1+s_2}}-\frac{1}{p^{1+s_3} }+\frac{1}{p^{1+s_1+s_2}}\right)\\ &=h_1(J,s_1,s_2,s_3)h_2(J,0) ,\end{split}\label{5.1}
\end{equation}
where $h_1$ is given by \eqref{3.12} and $h_2$ is given by \eqref{3.13}. On factoring out the dominant zeta factors from the product above we write
\begin{equation} F(s_1,s_2,s_3) = \frac{\zeta(1+s_1+s_2 ) }{\zeta(1+s_1) \zeta(1+s_2)\zeta(1+s_3)}h_1(J,s_1,s_2,s_3)h_j(s_1,s_2,s_3) \label{5.2}\end{equation}
which by \eqref{5.1} serves to define $h_j(s_1,s_2,s_3)$. By the same argument used for $h(s_1,s_2)$ in \eqref{2.15} we see that $h_j(s_1,s_2,s_3)$ converges absolutely and is analytic in the half-planes $\sigma_n > -c$, $1\le n\le 3$, for $c$ a positive constant ($c=\frac{1}{4}$ is acceptable). We also see immediately that $h_1(J,0,0,0)h_j(0,0,0) =  \gs(\text{\boldmath$j$})$, with $\text{\boldmath$j$}=(0,j)$ since we obtain the same product here as in the equation above \eqref{2.34}. Next, we define
\begin{equation} F(s_1,s_2,s_3) = k(s_1,s_2,s_3)\frac{s_1 s_2}{(s_1+s_2)\zeta(1+s_3)} ,\label{5.3} \end{equation}
and let \begin{equation}  \mathcal{J}(c_1,c_2,c_3)  = \frac {1}{ (2\pi i)^3}\mathop{\int}_{(c_3)\ } \!\mathop{\int}_{(c_2)\ } \!\mathop{\int}_{(c_1)\ } \frac{k(s_1,s_2,s_3)R^{s_1+s_2}}{s_1s_2(s_1+s_2)}\left( \frac{R^{s_3}}{{s_3}^2\zeta(1+s_3)}  \right) \, ds_1\,ds_2\, ds_3\label{5.4}
\end{equation}
so that  by \eqref{3.8} for $c_3>c_2>c_1>0$
\begin{equation} T_3(j) = \mathcal{J}(c_1,c_2,c_3) .\label{5.5}\end{equation}
Now, in the region \eqref{4.6} for $C$ appropriately small
\begin{equation} h_1(J,s_1,s_2,s_3)h_j(s_1,s_2,s_3) \ll \prod_{p|j}\left(1 +O(\frac{1}{p^{1-\delta}})\right) \ll \exp(D(\log |2j|)^\delta),\label{5.6}\end{equation}
by \eqref{2.36},
and hence in this region provided in addition $|s_i|\gg 1$, $1\le i\le 3$, and $|s_1+s_2|\gg \frac{1}{\log(|s_1|+2)}$ we see
\begin{equation} F(s_1,s_2,s_3) \ll \exp(D(\log |2j|)^\delta) \log^2(|t_1|+2)\log(|t_2|+2)
\log(|t_3|+2) \label{5.7}\end{equation}

We move the contour $(c_3)$ to the left to $\mathcal{L}_3$ which crosses the simple pole at $s_3=0$ and gives
\begin{equation} T_3(j)  = \frac {1}{ (2\pi i)^2} \mathop{\int}_{(c_2)\ } \!\mathop{\int}_{(c_1)\ } \frac{k(s_1,s_2,0)R^{s_1+s_2}}{s_1s_2(s_1+s_2)} \, ds_1\,ds_2 +\mathcal{J}(c_1,c_2,\mathcal{L}_3).\label{5.8}\end{equation}
The first integral may now be handled identically with how we handled the corresponding integral in \eqref{2.22} except the analytic factor $k(s_1,s_2,0)$ has the value
\[  k(s_1,0,0)=k(0,s_2,0) = \gs(\text{\boldmath$j$}).\]
We thus obtain
\begin{equation} T_3(j)  = \gs(\text{\boldmath$j$})N \log R +\mathcal{D}_0(j) N +\mathcal{J}(c_1,c_2,\mathcal{L}_3),\label{5.9}\end{equation}
where $\mathcal{D}_0(j)$ also satisfies the bound in \eqref{5.6}. For the last term we proceed as we did following \eqref{2.24}. We arrange in order the contours $\mathcal{L}_1$ and $\mathcal{L}_2$ to the right of
$\mathcal{L}_3$. Moving first $(c_1)$ to $\mathcal{L}_1$ we encounter a simple pole at $s_1=0$, and next  movimg $(c_2)$ to the left to $\mathcal{L}_2$ we obtain residues at $s_2=0$ and $s_2=-s_1$. By \eqref{2.27},
\[ \mathop{\int}_{(\mathcal{L}_3)\ }\left| \frac{R^{s_3}}{{s_3}^2\zeta(1+s_3)}\right|\, dt_3 \ll  e^{-c'\sqrt{\log R}} \]
and therefore by \eqref{5.7} these terms satisfy the bound $\ll N e^{-c'\sqrt{\log R}}$ for $\log |2j| \ll \log R$.
The theorem now follows in this case from the estimate
\[\mathcal{J}(\mathcal{L}_1,\mathcal{L}_2,\mathcal{L}_3) \ll  N e^{-c'\sqrt{\log R}}, \]
which also follows from the bound in \eqref{5.7} and \eqref{2.27}.

We finally turn to the case when $\text{\boldmath$a$}=(1,1,1)$.
Thus we assume $j_1\neq j_2 \neq j_3$, and in this case the only infinite product in \eqref{3.11} is $h_3(\Delta_3)$. Removing the dominant zeta-factors we have
\begin{equation} h_3(\Delta_3) = \frac{1}{\zeta(1+s_1)\zeta(1+s_2)\zeta(1+s_3)}h_4(\Delta_3) ,\label{5.10}\end{equation}
where
\begin{equation} \begin{split} h_4(\delta) &= \prod_{p \ndiv \delta}\left( \frac{ 1 - \frac{1}{p^{1+s_1}} - \frac{1}{p^{1+s_2}}- \frac{1}{p^{1+s_3}}}{\big( 1 - \frac{1}{p^{1+s_1}}\big)\big( 1 - \frac{1}{p^{1+s_2}}\big)\big( 1 - \frac{1}{p^{1+s_3}}\big)}\right)\\ & \qquad \times \prod_{p | \delta}\left( 1 - \frac{1}{p^{1+s_1}} \right)^{-1}\left( 1 - \frac{1}{p^{1+s_2}} \right)^{-1}\left( 1 - \frac{1}{p^{1+s_3}} \right)^{-1} .\end{split} \label{5.11}\end{equation}
We conclude that
\begin{equation} F(s_1,s_2,s_3) = \frac{1}{\zeta(1+s_1)\zeta(1+s_2)\zeta(1+s_3)} h_{\text{\boldmath{$j$}}}(s_1,s_2,s_3) , \label{5.12}\end{equation}
where
\begin{equation} h_{\text{\boldmath{$j$}}}(s_1,s_2,s_3) =
h_1(\kappa)h_2(\kappa, J_{12})h_2(\kappa, J_{13})h_2(\kappa, J_{23})h_4(\Delta_3)\label{5.13}\end{equation}
is analytic in the region \eqref{4.6}. We evaluate $T_3$ in \eqref{3.8} by moving the three contours to the left edge of this region. The only singularities in this region are simple poles at $s_1=0$, $s_2=0$, and $s_3=0$, and therefore
we obtain
\begin{equation}\begin{split}
T_3(j_1,j_2,j_3) &= h_{\text{\boldmath{$j$}}}(0,0,0)  + \frac {1}{ (2\pi i)^3}\mathop{\int}_{\mathcal{L}\ }\! \mathop{\int}_{\mathcal{L}\ }\! \mathop{\int}_{\mathcal{L}\ }F(s_1,s_2,s_3)\frac{R^{s_1+s_2+s_3}}{{s_1}^2 {s_2}^2 {s_3}^2} \, ds_1\,
ds_2 \, ds_3 \\ & = h_{\text{\boldmath{$j$}}}(0,0,0)+O(\exp(D (\log 2 ||\boldsymbol{j}||)^\delta)) e^{-c'\sqrt{\log R}}),\end{split}\label{5.14}\end{equation}
where the error term is obtained from the estimate \eqref{2.19} together with the estimate in \eqref{5.6} which applies for $h_{\text{\boldmath$j$}}(s_1,s_2,s_3) $, and \eqref{2.27}.
Since
\begin{equation}\begin{split} h_{\text{\boldmath{$j$}}}(0,0,0) &= \prod_{p| \kappa}\left( 1 - \frac{1}{p}\right)\prod_{\substack{p\ndiv \kappa\\ p|\Delta_3}}\left( 1 - \frac{2}{p}\right)\prod_{p\ndiv \Delta_3}\bigg( \frac{1 - \frac{3}{p}}{\big(1 - \frac{1}{p}\big)^3}\bigg)\prod_{p| \Delta_3}\left( 1 - \frac{1}{p}\right)^{-3} \\ &
=\prod_p\left( 1 - \frac{1}{p}\right)^{-3}\left(1 - \frac{\nu_{\text{\boldmath{$j$}}}(p)}{p}\right),\end{split} \label{5.15}\end{equation}
where
\begin{equation*}
    \nu_{\text{\boldmath{$j$}}}(p) =
        \left\{
        \begin{array}{ll}
        1, &\text{if $p|\kappa$,}\\
                  2, & \text{if $p|\Delta_3 $, $p \not | \kappa$,} \\
              3, &\text{if $p\not | \Delta_3$},
        \end{array}
        \right.
    \end{equation*}
we see that $\nu_{\text{\boldmath{$j$}}}(p)$ agrees with the definition of $\nu_p(\text{\boldmath{$j$}})$. We conclude by \eqref{3.6} that, for $\log 2||\boldsymbol{j}|| \ll \log R$,
\begin{equation}
\mathcal{S}_3(j_1,j_2,j_3) = \gs(\text{\boldmath{$j$}}) +O(e^{-c'\sqrt{\log R}}),\label{5.16}\end{equation}
which completes the proof of Theorem \ref{Theorem3.1}.

\section{The $k$-correlations }
We start by considering the sum
\begin{equation} \mathcal{S}(\text{\boldmath$j$}) = \sum_{n=1}^N \Lambda_R(n+j_1)\Lambda_R(n+j_2)\cdots \Lambda_R(n+j_k),\label{6.1}\end{equation}
where $\text{\boldmath$j$} = (j_1,j_2,\ldots ,j_k)$ and the $j_i$'s are not necessarily distinct. Proceeding as before we find as in \eqref{2.6} and \eqref{3.6} that

\begin{equation}\begin{split}
\mathcal{S}(\text{\boldmath$j$})&= N \sum_{\substack{d_1,d_2,\ldots, d_k \le R\\ (d_r,d_s)| j_s-j_r, \, 1\le r<s\le k }}\frac{\prod_{i=1}^k\mu(d_i)\log\frac{R}{d_i}}{[d_1, d_2, \ldots , d_k ]} + O(R^k) \\ &
=N T_k(\text{\boldmath$j$}) + O(R^k). \end{split}\label{6.2} \end{equation}
We next decompose the $d_i$'s into relatively prime factors. Let $\mathcal{P}(k)$ be the set of all non-empty subsets of the set of $k$ elements $\{ 1, 2, \ldots , k\}$ (This is just the power set with the empty set removed.) For $\mathcal{B} \in \mathcal{P}(k)$, we let $\mathcal{P}_\mathcal{B}(k)$ denote the set of all members of $\mathcal{P}(k)$ for which $\mathcal{B}$ is a subset. Thus for example
if $k= 4$ then
\[\mathcal{P}_{\{1,2\}}(4) = \{ \{1,2\}, \{1,2,3\}, \{1,2,4\}, \{1,2,3,4\} \} \]
Since the $d_i$'s are squarefree we now decompose them into the relatively prime factors
\begin{equation}  d_i = \prod_{m \in \mathcal{P}_{\{i\}}(k)} b_m ,
\qquad 1\le i \le k ,\label{6.3}\end{equation}
where $b_m$ is the product of all the primes that precisely divide all the $d_i$'s for which $i \in m$, and none of the other $d$'s.  This decomposition is unique and the $2^k-1$ $b_m$'s are pairwise relatively prime with each other (Here we subscript with the set $m$ rather than a list of its elements, so that for example $b_{12}$ in the earlier sections now is written $b_{\{1,2\}}$). Letting $\mathcal{D}(\text{\boldmath$j$})$ denote the divisibility conditions
\begin{equation}  \prod_{m\in \mathcal{P}_{ \{s,t\}}(k)}b_m \ \Big|\  j_t-j_s, \qquad 1\le s<t\le k , \label{6.4}\end{equation}
we see on  substituting \eqref{6.3} into \eqref{6.2} and applying \eqref{2.8} that, letting $\# m$ denote the number of elements of $m$,
\begin{equation}\begin{split}  T_k(\text{\boldmath$j$}) &= \sum_{\substack{d_1,d_2,\ldots, d_k \le R\\ \mathcal{D}(\text{\boldmath$j$})}}\left(\prod_{m\in \mathcal{P}(k)}\frac{{\mu(b_m)}^{\# m}}{b_m}\right)\left(\prod_{i=1}^k  \log\frac{R}{d_i} \right)\\
&
= \frac {1}{ (2\pi i)^k}\mathop{\int }_{(c_k)\ }\! \cdots \!\mathop{\int}_{(c_1)\ } F(s_1,\ldots, s_k)\prod_{i=1}^k\frac{R^{s_i}}{ {s_i}^2} ds_i ,\label{6.5}\end{split}\end{equation}
where
\begin{equation} F(s_1, \ldots , s_k) = \sumprime_{\substack{b_m,\  m\in  \mathcal{P}(k)\\ \mathcal{D}(\text{\boldmath$j$})  }}\prod_{m\in \mathcal{P}(k)}\frac{{\mu(b_m)}^{\# m}}{b_m^{1+ \tau_m}}, \label{6.6} \end{equation}
and
\begin{equation} \tau_m =\sum_{i\in m}s_i.\label{6.7} \end{equation}
None of the $b_m$ with $m$ a singleton set are constrained by the divisibility conditions in \eqref{6.4}, and  other $b_m$'s may also not be constrained when $j_s = j_t$ and $s\neq t$. Let $\mathcal{Q}_{\text{\boldmath$j$}}(k)$ denote the collection of $m \in \mathcal{P}(k)$ for which $b_m$ is not constrained by the divisibility conditions in \eqref{6.4}; we will specify $\mathcal{Q}_{\text{\boldmath$j$}}(k)$ precisely below. We now write $F(s_1, \ldots , s_k)$ as an Euler product
\begin{equation}  F(s_1, \ldots , s_k) = \prod_{p}\Big( 1 +\sum_{m\in \mathcal{Q}_{\text{\boldmath$j$}}(k)}\frac{(-1)^{\# m}}{p^{1+\tau_m}} + f_{\text{\boldmath$j$}}(p;s_1,s_2,\ldots ,s_k)\Big).
\label{6.8} \end{equation}
The local factor $f_{\text{\boldmath$j$}}(p;s_1,s_2,\ldots ,s_k)=0$ if $p > \max |j_s-j_t|$, so we factor out the zeta factors corresponding to the unconstrained variables and write
\begin{equation}
F(s_1,\ldots, s_k) =\left( \prod_{m\in \mathcal{Q}_{\text{\boldmath$j$}}(k)} \zeta(1+\tau_m)^{(-1)^{\# m}}\right)h_{\text{\boldmath$j$}}(s_1,\ldots , s_k) ,\label{6.9}\end{equation}
where we see $ h_{\text{\boldmath$j$}}(s_1,s_2,\ldots s_k)$ is analytic in the region
\begin{equation} \sigma_j \ge -\frac{C}{\log(|t_j|+2)}, \quad 1\le j\le k. \label{6.10} \end{equation}

Our analysis is complicated by the different cases of non-distinctness in  the numbers $j_1 , j_2, \ldots, j_k$. Suppose there are $r$ distinct values among the $j_i$'s. By reordering we may take these $r$ distinct values to be $j_1, j_2, \ldots, j_r$ which occur with multiplicity $a_1,a_2, \ldots, a_r$. We then have among the $k$ original $j$'s the relations, for $1\le i\le r$,
\begin{equation} j_i = j_{s_2(i)} = j_{s_3(i)} = \cdots =j_{s_{a_i}(i)}, \qquad r<s_2(i)<s_3(i) < \cdots < s_{a_i(i)}\le k. \label{6.11} \end{equation}
We also for notational purposes let
\begin{equation} s_1(i) = i , \qquad 1\le i \le r.\label{6.12}\end{equation}
Our divisibility conditions $\mathcal{D}(\text{\boldmath$j$})$ now reduce to the divisibility of the numbers $j_u-j_v$ when  $1\le u<v\le r$, which will be divisible by all the $b_m$'s with $m$ any element of the set
\begin{equation} {\mathcal{P}'}_{\{u,v\}}(k) = \bigcup_{n=1}^{a_v}\bigcup_{m=1}^{a_u}\mathcal{P}_{\{s_m(u),s_n(v)\}}(k). \label{6.13} \end{equation}
Hence $\mathcal{D}(\text{\boldmath$j$})$ is given by the equations
\begin{equation}  \prod_{m\in {\mathcal{P}'}_{ \{u,v\}}(k)}b_m \ \Big|\  j_v-j_u, \, \qquad 1\le u<v\le r . \label{6.14}\end{equation}
The indexing set $\mathcal{Q}_{\text{\boldmath$j$}}(k)$ for the unconstrained variables is now given by
\begin{equation} \mathcal{Q}_{\text{\boldmath$j$}}(k) = \bigcup_{i=1}^r \mathcal{Q}_i, \quad \mathrm{where} \ \mathcal{Q}_i = \left\{ m \neq \emptyset : m \subset \{ s_1(i),s_2(i),s_3(i), \ldots , s_{a_i}(i)\}\right\}, \label{6.15} \end{equation}
where we see the $\mathcal{Q}_i$ are disjoint from each other.
As we saw in the cases of $k=2$ and $k=3$ we will approximate $F(s_1, \ldots , s_k)$ in \eqref{6.9} by replacing each zeta factor with its first Taylor or Laurent term except when $m$ is a singleton that does not occur as a subset of any other $m'\in \mathcal{Q}_{\text{\boldmath$j$}}(k)$, in which case we retain the zeta factor. This latter situation occurs only when $a_i=1$ and in this case $\mathcal{Q}_i=\{ \{ i \}\}$. By reordering we can assume the multiplicities $a_i$ form a non-decreasing sequence, and we define $\omega$  to be the number of $a_i$'s for which $a_i =1$, and take $\omega =0$ if there are no such $a_i$'s. Thus
\begin{equation} a_1=a_2=\cdots =a_\omega =1, \quad  2 \le a_{\omega+1} \le a_{\omega +2} \le \ldots \le a_r .\label{6.16} \end{equation}

With this preparation, we now define $k(s_1,s_2,\ldots , s_k)$ by
\begin{equation}  F(s_1, \ldots , s_k) = k(s_1,\ldots ,s_k)\Big(\prod_{i=1}^\omega\frac{1}{\zeta(1+s_i)}\Big)\prod_{\substack{m\in \mathcal{Q}_{\text{\boldmath$j$}}(k)\\ m\neq \{i\} , 1\le i \le \omega}}(\tau_m)^{(-1)^{\# m+1}} .\label{6.17}\end{equation}
We thus see $k(s_1, \ldots , s_k)$ is analytic in the region \eqref{6.10} except possibly when $\tau_m =0$, and therefore by \eqref{6.5}
\begin{equation}\begin{split}  T_k(\text{\boldmath$j$}) =  \frac {1}{ (2\pi i)^k}\mathop{\int}_{(c_k)\ }\!\cdots \mathop{\int}_{(c_2)\ } \!\mathop{\int}_{(c_1)\ } k(s_1,\ldots ,s_k) &\Bigg(\prod_{\substack{m\in \mathcal{Q}_{\text{\boldmath$j$}}(k)\\ m\neq \{i\} , 1\le i \le \omega}}(\tau_m)^{(-1)^{\# m+1}}\Bigg)\\&
\qquad \times \Big(\prod_{i=1}^\omega\frac{R^{s_i}\, ds_i}{{s_i}^2\zeta(1+s_i)}\Big)\left(\prod_{i=\omega +1}^k\frac{R^{s_i}}{ {s_i}^2} ds_i\right). \end{split}\label{6.18}
\end{equation}
As before we see that in the region \eqref{6.10} if $|\tau_m|\gg \prod_{i\in m}\frac{1}{\log(|s_i|+2)}$ for $m \in \mathcal{P}(k)$ then
\begin{equation} F(s_1,s_2,\ldots , s_k) \ll_k \exp(D(\log 2||\boldsymbol{j}||)^\delta) \prod_{1\le i\le k} \log^B(|s_i|+2)\label{6.19}\end{equation}
where the $\text{\boldmath$j$}$ dependence may be estimated as in \eqref{5.6}.
We first move the contours $(c_1)$, $(c_2)$, \dots , $(c_\omega)$ to $\mathcal{L}$ where we encounter only simple poles at $s_i=0$. For the remaining contours we proceed to move each $(c_j)$ to $\mathcal{L}_j$ and encounter various poles which generate new terms which need to be integrated over the remaining contours. In these new terms we continue to move each $(c_j)$ either to the left to $\mathcal{L}_j$ or to the right to the reflected contour ${\mathcal{L}_j}^*$ according to whether the factor $R^{as_j}$ has $a>0$ or $a<0$.  This process eventually leads to a series of residues and integrals which are over contours  $\mathcal{L}_j$ or ${\mathcal{L}_j}^*$. These integrals are all estimated as before. We  arrange the contours so that they are well spaced where the conditions for \eqref{6.19} are satisfied, and apply \eqref{2.27} to see these integrals are, for $\log ||2\boldsymbol{j}|| \ll \log R$,
\[  \ll_k e^{-c'\sqrt{\log R}}. \]
We conclude that
\[ T_k(\text{\boldmath$j$}) = \sum_{i=0}^{\nu_{\text{\boldmath$j$}}(k)}\mathcal{D}_i(\text{\boldmath$j$})\log^iR +O_k( e^{-c'\sqrt{\log R}}).\]
Here  $\mathcal{D}_i(\text{\boldmath$j$})\ll_k \exp(D(\log 2||\boldsymbol{j}||)^\delta)$ since these terms are obtained from derivatives of $h_{\text{\boldmath$j$}}$ which will satisfy the same type of bound as in \eqref{5.6}. The leading residue terms will be those terms that contain the factor $k(0,0,\ldots , 0)$, and
hence are obtained by replacing the factor $k(s_1,s_2, \ldots , s_k)$ by  $h_{\text{\boldmath$j$}}(0,\ldots , 0)$ in the integral in \eqref{6.18}.  Because the $\mathcal{Q}_i$'s are disjoint, the new integrand without $k(s_1,s_2,\ldots,s_k)$ may be broken into a product of terms with disjoint variables which may be integrated separately, and thus the main term is
\begin{equation}\begin{split}
 h_{\text{\boldmath$j$}}(0,\ldots , 0)\prod_{i=\omega+1}^r & \left(\frac{1}{ (2\pi i)^{a_i}}\mathop{\int}_{(c_{a_i})\ }\!\cdots \mathop{\int}_{(c_{1})\ }
\Big( \prod_{m\in \mathcal{Q}_{i}}(\tau_m)^{(-1)^{\# m+1}}\Big) \prod_{i=1}^{a_i}\frac{R^{s_i}}{ {s_i}^2} ds_i \right) \\& \qquad  \times \left( \frac{1}{ 2\pi i} \mathop{\int}_{(c)\ }\!  \frac{R^s}{s^2\zeta(1+s)} ds\right)^\omega .\label{6.20} \end{split}\end{equation}
On moving $(c)$ to $\mathcal{L}$ we see by \eqref{2.27} that
\[ \frac{1}{ 2\pi i} \mathop{\int}_{(c)\ }\!  \frac{R^s}{s^2\zeta(1+s)} ds = 1 +O(e^{-c'\sqrt{\log R}}),\]
 and  therefore these factors may be ignored in \eqref{6.20} with an error of $O(e^{-c'\sqrt{\log R}})$.  As before we make a change of variable to remove the $R$ dependence in the first integral but retain the $s_j$'s as our variables in the $\tau_m$'s through the replacement $s_j \rightarrow s_j/\log R $ and $(c_j)\rightarrow (\frac{c_j}{\log R})$. We thus have the main term above is
\begin{equation} \mathcal{C}_k(\text{\boldmath$j$}) h_{\text{\boldmath$j$}}(0,\ldots , 0) \log^{\nu_{\text{\boldmath$j$}}(k)} R , \label{6.21}
\end{equation}
where
\begin{equation} \mathcal{C}_k(\text{\boldmath$j$})= \prod_{i=\omega +1}^r \mathcal{C}_{a_i},  \label{6.22} \end{equation}
\begin{equation} \mathcal{C}_k = \mathcal{C}_k(k) = \frac {1}{ (2\pi i)^{k}}\mathop{\int}_{(c_k)\ }\!\cdots \!\mathop{\int}_{(c_{1})\ }\Big( \prod_{m \in \mathcal{P}(k)}(\tau_m)^{(-1)^{\# m+1}}\Big) \prod_{i=1}^k\frac{e^{s_i}}{ {s_i}^2} ds_i .\label{6.23}
\end{equation}
The power $\nu_{\text{\boldmath$j$}}(k)$ of $\log R$ obtained by the change of variables  is equal to the number of $\tau_m$ in the denominator of the integrand  minus the number of $\tau_m$ in the numerator plus $k -\omega$ from the product of $R^{s_i}/{s_i}^2$. Thus
\begin{equation}\begin{split} \nu_{\text{\boldmath$j$}}(k)&= k - \omega + \sum_{\substack{m\in \mathcal{Q}_{\text{\boldmath$j$}}(k)\\ m\neq \{i\}, 1\le i\le \omega}}(-1)^{\# m } \\&
=  k  + \sum_{m\in \mathcal{Q}_{\text{\boldmath$j$}}(k)}(-1)^{\# m }.\end{split} \label{6.24}
\end{equation}

To evaluate the sum above, we use \eqref{6.15} and count the elements of $\mathcal{Q}_{\text{\boldmath$j$}}(k)$ according to their cardinality. Thus
\begin{equation}\begin{split} \sum_{m\in \mathcal{Q}_{\text{\boldmath$j$}}(k)}(-1)^{\# m }
&=  \sum_{i=1}^r\left( - \binom{a_i}{1} + \binom{a_i}{2}-\binom{a_i}{3} + \cdots + (-1)^i\binom{a_i}{a_i} \right) \\
&= - \sum_{i=1}^r 1 \\
&=-r ,\end{split} \label{6.25}\end{equation}
and therefore we have
\begin{equation}  \nu_{\text{\boldmath$j$}}(k)= k-r .\label{6.26}\end{equation}
Next, from \eqref{6.8} and \eqref{6.9} we have, using \eqref{6.25},
\[ \begin{split} h_{\text{\boldmath$j$}}(0,0,\ldots ,0) &= \prod_{p}\left(1 - \frac{1}{p}\right)^{\displaystyle{\sum_{m\in \mathcal{Q}_{\text{\boldmath$j$}}(k)}(-1)^{\# m }}}\left(1 + \frac{1}{p}\sum_{m\in \mathcal{Q}_{\text{\boldmath$j$}}(k)}(-1)^{\# m }+f_{\text{\boldmath$j$}}(p;0,0,\ldots ,0)
\right)\\ & = \prod_{p}\left(1-\frac{1}{p}\right)^{-r}\left(1 - \frac{r}{p} + f_{\text{\boldmath$j$}}(p;0,0,\ldots ,0)
\right).\end{split}\]
Letting
\begin{equation} \text{\boldmath$j'$} = (j_1,j_2,\ldots j_r), \label{6.27}\end{equation}
we will now prove that
\begin{equation} f_{\text{\boldmath$j$}}(p;0,0,\ldots ,0)= \frac{r-\nu_p(\text{\boldmath$j'$})}{p} \label{6.28}\end{equation}
and thus by \eqref{1.5} we conclude
\begin{equation} h_{\text{\boldmath$j$}}(0,0,\ldots ,0) = \gs(\text{\boldmath$j'$}).\label{6.29}\end{equation}
From \eqref{6.6}, \eqref{6.8}, and \eqref{6.14} we see that
\[ f_{\text{\boldmath$j$}}(p;s_1,s_2,\ldots ,s_k)=\sum_{1\le u<v\le r}\sum_{\substack{m\in {\mathcal{P}'}_{ \{u,v\}}(k)\\ p|j_v -j_u}} \frac{(-1)^{\# m}}{p^{1+\tau_m}},\]
and hence
\[f_{\text{\boldmath$j$}}(p;0,\ldots ,0)=\frac{1}{p}\sum_{1\le u<v\le r}\sum_{\substack{m\in {\mathcal{P}'}_{ \{u,v\}}(k)\\ p|j_v -j_u}} (-1)^{\# m}.\]
Thus, the proof of \eqref{6.28} reduces to proving
\[ \sum_{1\le u<v\le r}\sum_{\substack{m\in {\mathcal{P}'}_{ \{u,v\}}(k)\\ p|j_v -j_u}} (-1)^{\# m}= r - \nu_p(\text{\boldmath$j'$}) .\]
If $\nu_p(\text{\boldmath$j'$})=q$, then  the $r$ distinct numbers $j_1,j_2, \ldots , j_r$ must fall into $q$ residue classes congruent to numbers $r_1, r_2, \ldots, r_q $ modulo $p$, say with multiplicities $m_1,m_2, \ldots , m_q$, with $m_i\ge 1$ and $m_1+m_2 +\cdots +m_q =r$. Let
\[ \mathcal{M}_i = \{ n: j_n\equiv r_i  (\mathrm{mod}\ p), \ 1\le n\le r\}, \]
so that $\#\mathcal{M}_i =m_i$. Next, recalling \eqref{6.15}, let
\[\begin{split} \tilde{\mathcal{M}}_i &= \{ s_1(n), s_2(n), \ldots , s_{a_n}(n) : n \in \mathcal{M}_i \} \\ & =
\bigcup_{n\in \mathcal{M}_i}\bigcup_{w=1}^{a_n} s_w(n), \end{split}\]
and
\[ \tilde{\mathcal{P}}(i) = \{ m \neq \emptyset : m \subset  \tilde{\mathcal{M}}_i , \#m \ge 2\},\]
so that
\[  \# \tilde{\mathcal{M}}_i  = \sum_{n\in \mathcal{M}_i}a_n.\]
 The condition $p|j_v-j_u$ will only hold if $j_v$ and $j_u$ are in the same residue class, and hence if and only if
\[ \{u,v\} \in \bigcup_{i=1}^q  \tilde{\mathcal{P}}(i).\]
 Hence we conclude
\[ \bigcup_{\substack { 1\le u<v\le r\\ p|j_v-j_u}}{\mathcal{P}'}_{ \{u,v\}}(k) = \bigcup_{i=1}^q \tilde{\mathcal{P}}(i)  -  \Big( \mathcal{Q}_{\text{\boldmath$j$}}(k) - \{\{1\},\{2\},\ldots , \{k\} \} \Big) \]
since  $\tilde{\mathcal{P}}(i)$ also contains all the unconstrained sets of variables with $\ge 2$ elements which are not included in ${\mathcal{P}'}_{ \{u,v\}}(k)$ and
$\mathcal{Q}_{\text{\boldmath$j$}}$ also includes the singleton sets which are not in $\tilde{\mathcal{P}}(i)$. We conclude by \eqref{6.25}
\[ \begin{split}\sum_{1\le u<v\le r}\sum_{\substack{m\in {\mathcal{P}'}_{ \{u,v\}}(k)\\ p|j_v -j_u}} (-1)^{\# m} &= \sum_{i=1}^q\sum_{m\in  \tilde{\mathcal{P}}(i)}(-1)^{\# m}
  - \sum_{m\in \mathcal{Q}_{\text{\boldmath$j$}}(k)}(-1)^{\# m} -\sum_{m=1}^k 1 \\& =\sum_{i=1}^q\sum_{j=2}^{\#\tilde{\mathcal{M}}_i}\binom{\#\tilde{\mathcal{M}}_i}{j} (-1)^j +r-k \\
&= \sum_{i=1}^q\Big( -1 + \#\tilde{\mathcal{M}}_i\Big)  +r-k \\
&= -q + \sum_{i=1}^q \sum_{n\in \mathcal{M}_i}a_n +r-k\\
&= r-q \\&
= r - \nu_p(\text{\boldmath$j'$}),\end{split}\]
which proves \eqref{6.28}.

Combining our results we have shown that
\begin{equation}T_k(\text{\boldmath$j$})
 = \mathcal{C}_k(\text{\boldmath$j$}) \gs(\text{\boldmath$j'$}) \log^{k-r}R  + \sum_{i=1}^{k-r}\mathcal{D}_i(\text{\boldmath$j$})\log^{k-r-i}R + O_k( e^{-c'\sqrt{\log R}}). \label{6.30} \end{equation}
Returning to the notation of \eqref{1.2}, we see in \eqref{6.2} that
  $\mathcal{S}(\text{\boldmath$j$})
=\mathcal{ S}_k(N,\text{\boldmath$j'$}, \text{\boldmath $a$})$ with $\text{\boldmath$j'$}$ given by \eqref{6.27}. Thus we have proved the following refined version of Theorem \ref{Theorem1}.
\begin{theorem} \label{Theorem6} Given $k\ge 1$,  let $\text{\boldmath$j$} = (j_1,j_2, \ldots , j_r)$ and $\text{\boldmath$a$} = (a_1,a_2, \ldots a_r)$, where the $j_i$'s are distinct integers, and $a_i\geq 1$ with $\sum_{i=1}^r a_i = k$. Assume $\log 2||\boldsymbol{j}|| \ll \log R$. Then
we have
\begin{equation}\begin{split} \mathcal{S}_k(N,\text{\boldmath$j$},\text{\boldmath$a$}) = \mathcal{ C}_k(\text{\boldmath$a$})\gs(\text{\boldmath$j$})N(\log R)^{k-r} &+ \sum_{i=1}^{k-r}\mathcal{D}_i(\text{\boldmath$j$})N \log^{k-r-i}R \\ & \qquad + O_k(||\text{\boldmath$j$}^\epsilon || Ne^{-c'\sqrt{\log R}})
+O(R^k),\label{6.31}\end{split} \end{equation}
where $\mathcal{D}_i(\text{\boldmath$j$})\ll \exp(D(\log 2||\boldsymbol{j}||)^\delta)$.
\end{theorem}

The constants $\mathcal{C}_k(\text{\boldmath$a$})$ satisfy the relation \eqref{1.7} because we may extend the product in \eqref{6.22} to include the terms with $1\le i\le \omega$ since for these terms $\mathcal{C}_{a_i}=\mathcal{C}_1 =1$.

\section{The constants $\mathcal{C}_k$ }

We now discuss the constants $\mathcal{C}_k$ that arise in the $k$-correlation results. We will show the integrals defining these constants converge absolutely and vanish as the contours are moved to the left to $-\infty$. New terms generated from the residues can be handled similarly by moving the contours appropriately to the left or right to infinity. Therefore the $\mathcal{C}_k$'s can be calculated by a residue calculation. We have already computed $\mathcal{C}_2$ and $\mathcal{C}_3$ before, but we will recompute them directly and also compute $\mathcal{C}_4$.

Consider first
\begin{equation} \mathcal{C}_2 = \mathcal{C}_2(2) = \frac {1}{ (2\pi i)^2}\mathop{\int }_{(c_2)\ }\! \mathop{\int}_{(c_1)\ } \frac{e^{s_1+s_2}}{s_1s_2(s_1+s_2)} \, ds_1\,ds_2 ,\label{7.1} \end{equation}
for $c_2>c_1>0$.
Letting
\begin{equation} I(\sigma_1,\sigma_2) = \frac{1}{(2\pi i)^2}\mathop{\int }_{(\sigma_2)\ }\! \mathop{\int}_{(\sigma_1)\ } m_2(s_1,s_2) \, ds_1\, ds_2 ,\label{7.2}\end{equation}
where
\begin{equation} m_2(s_1,s_2) =\frac{ e^{s_1+s_2}}{s_1s_2(s_1+s_2)} ,\label{7.3}\end{equation}
we may rewrite \eqref{7.1} as
\[  \mathcal{C}_2 =  I(\sigma_1, \sigma_2), \quad  \sigma_2>\sigma_1 >0.\]
If
\begin{equation} |s_1|\gg 1, \quad |s_2|\gg 1, \quad |s_1+s_2| \gg 1, \label{7.4} \end{equation}
then we have the bound
\begin{equation}  m_2(s_1,s_2) \ll \frac{e^{\sigma_1+\sigma_2}}{(1+|t_1|)(1+|t_2|)(1+|t_1+t_2|)}.\label{7.5} \end{equation}
Suppose the contours $(\sigma_1)$ and $(\sigma_2)$ are arranged so that the conditions in \eqref{7.4} hold; this can always be done in what follows. Then by symmetry we may estimate $I(\sigma_1,\sigma_2)$ over the restricted integration range $|s_1|\le |s_2|$ and $t_2 \ge 0$, and hence
\begin{equation} \begin{split} I(\sigma_1,\sigma_2) & \ll e^{\sigma_1+\sigma_2}\int_0^\infty\frac{1}{1+t_2} \left(\int_{-t_2}^{t_2} \frac{dt_1}{(1 + |t_1|)(1+t_1+t_2)}\right)\, dt_2 \\ & \ll e^{\sigma_1+\sigma_2}\int_0^\infty\frac{\log(1+t_2)}{(1+t_2)^2}\, dt_2\\ &\ll e^{\sigma_1+\sigma_2} .\end{split}\label{7.6} \end{equation}
From this estimate we see $\mathcal{C}_2$ is well-defined, and further that
\begin{equation}  \lim_{\sigma_j \to -\infty} I(\sigma_1,\sigma_2) = 0, \quad  j=1,2.\label{7.7}\end{equation}
Since $m_2(s_1,s_2) \to 0$ if $|t_1|\to \infty$ or $|t_2|\to \infty$, we can move the contour integrals defining $\mathcal{C}_2$  to the left and $\mathcal{C}_2$ can be evaluated from the residues encountered.

Moving first $c_1 \to - \infty $, we encounter residues at the simple poles $s_1=0$ and $s_1=-s_2$, and hence
\begin{equation}
\mathcal{C}_2= \frac {1}{ 2\pi i}\mathop{\int}_{(c_2) \ }\left( \frac{e^{s_2}}{{s_2}^2} - \frac{1}{{s_2}^2} \right)\, ds_2 .  \label{7.8} \end{equation}
Next moving $c_2 \to -\infty$, we encounter a pole at $s_2=0$ with residue equal to 1. Hence
\begin{equation} \mathcal{C}_2(2)=1. \label{7.9} \end{equation}

Turning next to $\mathcal{C}_3$, we have for $\sigma_3>\sigma_2>\sigma_1 >0 $
\begin{equation}  \mathcal{C}_3 = \mathcal{C}_3(3)= I(\sigma_1,\sigma_2,\sigma_3) \label{7.10}\end{equation}
where
\begin{equation} I(\sigma_1,\sigma_2,\sigma_3) := \frac {1}{ (2\pi i)^3}\mathop{\int}_{(\sigma_3)\ }\!\mathop{\int}_{(\sigma_2)\ } \!\mathop{\int}_{(\sigma_1)\ } m_3(s_1,s_2,s_3) \, ds_1\,ds_2\, ds_3 ,\label{7.11} \end{equation}
\begin{equation} m_3(s_1,s_2,s_3)= \frac{(s_1+s_2+s_3)e^{s_1+s_2+s_3}}{s_1s_2s_3(s_1+s_2)(s_1+s_3)(s_2+s_3)}. \label{7.12}\end{equation}
If
\begin{equation} |s_i|\gg 1, \ 1\le i \le 3,\quad |s_i+s_j|\gg 1, \ 1\le i<j\le 3, \label{7.13}\end{equation}
then
\[ \left|\frac{s_1+s_2+s_3}{s_2(s_1+s_3)}\right| = \left| \frac{1}{s_2}+ \frac{1}{s_1+s_3}\right| \ll 1,\]
 and hence
\begin{equation}  m_3(s_1,s_2,s_3)\ll e^{\sigma_1+\sigma_2+\sigma_3}\frac{1}{|s_1||s_3||s_1+s_2||s_2+s_3|} .\label{7.14}\end{equation}
By symmetry  it sufficies to estimate $I$ over the region $|s_1|\le |s_2|\le |s_3|$,
and hence as in \eqref{7.6}, assuming \eqref{7.13} holds,
\begin{equation}\begin{split}  I(\sigma_1,\sigma_2,\sigma_3) &\ll e^{\sigma_1+\sigma_2+\sigma_3} \mathop{\int}_{(\sigma_3)}\frac{dt_3}{|s_3|}   \mathop{\int}_{\substack{(\sigma_2)\\ |s_2|\le s_3}}\frac{1}{|s_2+s_3|}dt_2 \mathop{\int}_{\substack{(\sigma_1)\\ |s_1|\le |s_2|} } \frac{1}{|s_1||s_1+s_2|} \, dt_1 \\&
\ll e^{\sigma_1+\sigma_2+\sigma_3} \mathop{\int}_{(\sigma_3)}\frac{dt_3}{|s_3|}  \mathop{\int}_{\substack{(\sigma_2)\\ |s_2|\le s_3}}\frac{\log(2+|t_2|)}{|s_2||s_2+s_3|}dt_2 \\&
\ll e^{\sigma_1+\sigma_2+\sigma_3} \mathop{\int}_{(\sigma_3)}\frac{\log^2(2+|t_3|)}{|s_3|^2}dt_3
\\ & \ll e^{\sigma_1+\sigma_2+\sigma_3} .\label{7.15}\end{split}\end{equation}
Thus $\mathcal{C}_3$ is well defined and
\begin{equation}  \lim_{\sigma_j \to -\infty} I(\sigma_1,\sigma_2,\sigma_3) = 0, \quad  j=1,2,3.\label{7.16}\end{equation}
Since $m(s_1,s_2,s_3)$ vanishes as $|s_i| \to \infty$ for any $1\le i\le 3$, we can move our contours to the left and evaluate $\mathcal{C}_3$ from the residues encountered.

We now evaluate $\mathcal{C}_3$ by making a simplification suggested by D. Farmer which simplifies  the calculation of the residues in this case. Observe that if we break the integral defining $\mathcal{C}_3$  into three pieces corresponding to each term in the numerator, then by symmetry these terms are equal to each other, and therefore
\begin{equation} \mathcal{C}_3= \frac {3}{ (2\pi i)^3}\mathop{\int}_{(c_3)\ }\!\mathop{\int}_{(c_2)\ } \!\mathop{\int}_{(c_1)\ } \frac{e^{s_1+s_2+s_3}}{s_2s_3(s_1+s_2)(s_1+s_3)(s_2+s_3)} \, ds_1\,ds_2\, ds_3 .\label{7.17} \end{equation}
We choose $c_j=j$ for $1\le j\le 3$. Moving first $c_1 \to - \infty $, we encounter the simple poles  $s_1=-s_2$, and $s_1 = -s_3$ and obtain
\[ \mathcal{C}_3= \frac {3}{ (2\pi i)^2}\mathop{\int}_{(c_3)\ }\!\mathop{\int}_{(c_2)\ }  \frac{e^{s_2} - e^{s_3}}{s_2s_3(s_2+s_3)(s_2-s_3)}\,ds_2\, ds_3 .\]
As a function of $s_2$ the integrand has simple poles at $s_2=0$, $s_2=s_3$, and $s_2=-s_3$. On moving the contour $c_2 \to -\infty$ we do not encounter the pole at $s_2=s_3$, since $s_3$ is on $(c_3) =(3)$ and is to the right of $(c_2)=(2)$. The residue at $s_2 =0$ is
\[ \frac{e^{s_3}-1}{{s_3}^3} , \]
and the residue at $s_2=-s_3$ is
\[ \frac{e^{-s_3}-e^{s_3}}{2{s_3}^3} ; \]
therefore we have
\[ \mathcal{C}_3= \frac {3}{ 2\pi i}\bigg(\mathop{\int}_{(c_3)\ } \bigg( \frac{e^{s_3}}{2{s_3}^3}-\frac{1}{{s_3}^3}\bigg)\, ds_3 +\mathop{\int}_{(c_3)\ }  \frac{e^{-s_3}}{2{s_3}^3}\, ds_3\bigg).\]
In the first integral above we move $c_3 \to -\infty$ and encounter a triple pole at $s_3=0$ with residue equal to $\frac{1}{4}$. The second integral is equal to zero since we encounter no singularities when we move $c_3 \to \infty$ where the integral vanishes. We conclude
\begin{equation} \mathcal{C}_3=\frac{3}{4}. \label{7.18} \end{equation}

We now turn to the general case, and have for $\sigma_1, \sigma_2, \ldots ,\sigma_k>0$
\begin{equation} \mathcal{C}_k = I_k(\sigma_1, \sigma_2, \ldots ,\sigma_k) \label{7.19} \end{equation}
where
\begin{equation} I_k(\sigma_1, \sigma_2, \ldots ,\sigma_k)=\frac {1}{ (2\pi i)^{k}}\mathop{\int}_{(\sigma_k)\ }\!\cdots \!\mathop{\int}_{(\sigma_{1})\ }m_k(s_1,s_2,\ldots , s_k) ds_1\,ds_2\ldots ds_k, \label{7.20} \end{equation}
\begin{equation} m_k(s_1, s_2, \ldots s_k)=\Big( \prod_{m \in \mathcal{P}(k)}(\tau_m)^{(-1)^{\# m+1}}\Big) \prod_{i=1}^k\frac{e^{s_i}}{ {s_i}^2},\label{7.21} \end{equation}
and $\tau_m$ is defined in \eqref{6.7}. We will prove that $I_k$ is well defined and converges by proving that provided $|\tau_m |\gg 1$ for $m\in \mathcal{P}(k)$ then
\begin{equation}I_k(\sigma_1, \sigma_2, \ldots \sigma_k)\ll_k e^{\sigma_1+\sigma_2+\cdots +\sigma_k}. \label{7.22}\end{equation}
We have already proved this for $k=2$ and $k=3$ and we will proceed inductively. Define a function $g_k(s_1,s_2,\ldots , s_k)$ by first taking $g_2(s_1,s_2)=1$, and define $g_3(s_1,s_2,s_3)$ by
\[m_3(s_1,s_2,s_3)= \frac{ e^{s_1+s_2+s_3}}{s_1(s_1+s_2)(s_2+s_3)s_3}g_3(s_1,s_2,s_3) \]
so that
\[g_3(s_1,s_2,s_3) = \frac{s_1+s_2+s_3}{s_2(s_1+s_3)} .\]
We next define $g_4(s_1,s_2,s_3,s_4)$ by
\[ m_4(s_1,s_2,s_3,s_4) =\frac{ e^{s_1+s_2+s_3+s_4}}{s_1(s_1+s_2)(s_2+s_3)(s_3+s_4)s_4} g_3(s_1,s_2,s_3)g_4(s_1,s_2,s_3,s_4)\]
so that
\[  g_4(s_1,s_2,s_3,s_4)= \frac{ (s_1+s_2+s_4)(s_1+s_3+s_4)(s_2+s_3+s_4)}{s_3(s_1+s_4)(s_2+s_4)(s_1+s_2+s_3+s_4)}.\]
We define $g_k$ in general by
\begin{equation} m_k(s_1,s_2,\ldots , s_k) = \left(\prod_{i=1}^k e^{s_i}\right)\frac{1}{s_1s_k}\left(\prod_{i=1}^{k-1}\frac{1}{s_i+s_{i+1}}\right)\prod_{i=3}^kg_i(s_1,s_2,\ldots , s_i) .\label{7.23} \end{equation}
We will show below that, provided $|\tau_m |\gg 1$ for $m\in \mathcal{P}(k)$,
\begin{equation} \prod_{i=3}^k|g_i(s_1,s_2,\ldots , s_k)| \ll_k 1. \label{7.24}\end{equation}
Assuming this estimate for the moment, we see
\begin{equation}I_k(\sigma_1, \sigma_2, \ldots \sigma_k) \ll_k\left( \prod_{i=1}^k e^{\sigma_i}\right) \mathop{\int}_{(\sigma_k)\ }\!\cdots \!\mathop{\int}_{(\sigma_{1})\ }\frac{1}{|s_1||s_k|}\left(\prod_{i=1}^{k-1}\frac{1}{|s_i+s_{i+1}|}\right) ds_1\,ds_2\ldots ds_k .\label{7.25} \end{equation}
As in \eqref{7.6} and \eqref{7.15}, provided $|s_i|\gg 1$ and $|s_i+s_{i+1}|\gg 1$ we see that
\begin{equation} \mathop{\int}_{|s_i|\le |s_{i+1}|}\frac{\log^B(2 + |t_{i})}{|s_i||s_{i}+s_{i+1}|}dt_i \ll \frac{\log^{B+1}(2+|t_{i+1}|)}{|s_{i+1}|} .\label{7.26}\end{equation}
By symmetry we may estimate the integral in \eqref{7.25} over the region $|s_1|\le |s_2|\le \cdots \le |s_k|$, and hence we see on using the estimate \eqref{7.26} repeatedly that the integral converges, which proves \eqref{7.22}. It remains to prove \eqref{7.24}.

We have already seen that $g_3 \ll 1$ in the equation above \eqref{7.14}. For $g_4$ we have
\[ \begin{split} |g_4(s_1,s_2,s_3,s_4)| &= \frac{|s_1+s_3+s_4|}{|s_3||s_1+s_4|}\cdot \frac{|s_1+s_2+s_4||s_2+s_3+s_4|}{|s_2+s_4||s_1+s_2+s_3+s_4|}\\& =
 \left| \frac{1}{s_3} + \frac{1}{s_1+s_4}\right| \left|\frac{1 +\frac{s_1}{s_2+s_4}}{1 + \frac{s_1}{s_2+s_3+s_4}}\right| \\ &
\ll 1,\end{split} \]
since $\tau_m \gg 1$.
For the general case we claim
\begin{equation} \begin{split} g_k(s_1, s_2, \ldots , s_k) &=\frac{(s_1+s_{k-1}+s_k)}{s_{k-1} (s_1+s_k)} \tilde g_k \\ & \qquad \times  \left( \frac{(s_1+s_2+ \cdots + s_{k-2}+s_k)(s_2+s_3+\cdots +s_k)}{(s_2 +s_3 +\cdots +s_{k-2}+s_k)(s_1+s_2+\cdots +s_k)}\right)^{(-1)^k},\end{split}\label{7.27}\end{equation}
and $\tilde g_k $ is a product of terms of the form
\begin{equation} \left( \frac{s_1 +\eta_m}{\eta_m}\cdot \frac{\eta_m +s_j}{s_1 +\eta_m + s_j} \right)^{(-1)^{\# m}}, \qquad  \eta_m = \sum_{i\in m}s_i, \label{7.28} \end{equation}
where $m\in \mathcal{P}(k)$, $m$ has between two and $k-3$ terms,  and $k\in m$ but $1\notin m$ since we have separated out these terms above.  For each $m$ we  may choose  $j$ arbitrarily except for the requirements that $2\le j\le k$ and this factor has not already appeared in the product. For example,
\[ \tilde g_5 = \frac{(s_1+s_2+s_5)(s_2+s_4+s_5)}{(s_2+s_5)(s_1+s_2+s_4+s_5)}\cdot \frac{(s_1+s_3+s_5)(s_3+s_4+s_5)}{(s_3+s_5)(s_1+s_3+s_4+s_5)} \]
For $k=6$ we can write $\tilde g_6$ as a product of 6 terms  corresponding to the pairs $(\eta_m ,j)$ given by $(s_2+s_6 , 3)$, $(s_3+s_6 , 4)$, $((s_4+s_6,2)$, $(s_3+s_5+s_6, 2)$, $(s_2+s_5+s_6, 4)$, and $(s_4+s_5+s_6, 3)$. In the general case, we need to check that we will run through all the $m \in \mathcal{P}(k)$ with this decomposition. Then, since
\begin{equation} \begin{split} \left( \frac{s_1 +\eta_m}{\eta_m}\cdot \frac{\eta_m +s_j}{s_1 +\eta_m + s_j} \right)^{(-1)^{\# m}}&= \left( \frac{1+ \frac{s_1}{\eta_m}}{1 + \frac{s_1}{\eta_m +s_j}}\right)^{(-1)^{\# m}}\\&
\ll_k 1 \label{7.29} \end{split}\end{equation}
if $|\tau_m|\gg 1$, and similarly for  the first and last individual terms in \eqref{7.27}, we see that \eqref{7.24} follows.
To establish that $\tilde g_k$ can be written as a product of terms of the form \eqref{7.28}, we use a counting argument to verify that in constructing the product we use all of the remaining terms $\tau_m$ with $m\in \mathcal{P}(k)$ in $m_k(s_1,s_2,\ldots , s_k)$.  We examine this by grouping according to $\# m$.  All of singleton elements  have already appeared in \eqref{7.23}. Note every $m$ for which $\eta_m$ occurs in  $\tilde g_k$ must contain $k$, because all the terms $\tau_m$ without $s_k$ have already appeared in some $g_{i}$, $2\le i\le k-1$, or the product in \eqref{7.23}. Further, $1$ can not be in $m$ because we have separated out terms with $s_1$ in \eqref{7.28}. Therefore, for $\eta_m$ with $\# m =2$, we can take $\eta_m = s_i +s_k$ with  $2\le i\le k-2$ since $s_{k-1} +s_{k}$ already appears  in the product in \eqref{7.23},  and hence there are  $k-3$ such terms. All of $\# m =2$ terms are now accounted for in \eqref{7.23}. For $\eta_m$ with $\# m =3$, we can take $\binom{k-2}{2}$ choices for $\eta_m$ in the first denominator of \eqref{7.28} except we can not use any of the $k-3$ factors $\eta_m +s_j$ we just obtained when $\# m=2$, and hence there are
\[\binom{k-2}{2} - (k-3) = \binom{k-3}{2} \]
such terms in \eqref{7.28}, and this accounts for all the $\tau_m$ with $\# m =3$. In general we see that there will be $\binom{k-3}{j-1}$ terms $\eta_m$ with $\# m =j$ occurring in \eqref{7.28} for $1\le j \le k-3$. Since $\# m =k-3$ is the largest size set for which $\eta_m$ occurs in \eqref{7.28}, we only need to check that the terms with $\# m =k-2$, $\# m =k-1$, and $\# m =k$ are accounted for in the earlier terms and the last term in \eqref{7.27}. There is only one element with $\# m =k$, and it only occurs in the last term in \eqref{7.27}. For the $k-1$ elements  with $\#m =k-1$, $k-3$ of them occur in the product formed from elements with $\# m =k-3$ in \eqref{7.28}, and the remaining two terms are in the last term of \eqref{7.27}. Finally, for the $\binom{k-1}{2}$ elements with $\# m =k-2$, we have $\binom{k-3}{k-5}$ of them in the terms in \eqref{7.28} coming from terms $\eta_m$ with $\# m =k-4$ and $2(k-3)$ of them from $\eta_m$ with  $\# m = k-3$. Since
\[ \binom{k-1}{2} -\left(\binom{k-3}{k-5} +2(k-3)\right) =1 \]
we see this leaves exactly the one term needed for the last term in \eqref{7.27}.

We conclude this section by evaluating $\mathcal{C}_4$.  We have
\begin{equation} \mathcal{C}_4= \frac {1}{ (2\pi i)^4}\mathop{\int}_{(c_4)\ }\!\mathop{\int}_{(c_3)\ }\!\mathop{\int}_{(c_2)\ } \!\mathop{\int}_{(c_1)\ }m_4(s_1,s_2,s_3,s_4)  \, ds_1\,ds_2\, ds_3 \, ds_4, \label{7.30}\end{equation}
where $c_4>c_3>c_2>c_1>0$ and
\[\begin{split} &m_4(s_1,s_2,s_3,s_4)=\\& \frac{(s_1+s_2+s_3)(s_1+s_2+s_4)(s_1+s_3+s_4)(s_2+s_3+s_4)e^{s_1+s_2+s_3+s_4}}{s_1s_2s_3s_4(s_1+s_2)(s_1+s_3)(s_1+s_4)(s_2+s_3)(s_2+s_4)(s_3+s_4)(s_1+s_2 +s_3+s_4)}.\end{split}\]
Moving $c_1$ to $-\infty$ we encounter poles at $s_1=0$, $s_1=-s_2$, $s_1=-s_3$, $s_1=-s_4$, and $s_1= -s_2-s_3-s_4$. The residues at the poles $s_1=-s_2$, $s_1=-s_3$, $s_1=-s_4$ are symmetric with each other, but because the contours are not symmetric they do not contribute equally.

First, the residue at $s_1=0$ is equal to
\[ \frac{e^{ s_2 +  s_3 +  s_4}}{
    {s_2}^2{ s_3}^2{ s_4}^2}, \]
and thus by \eqref{2.8} or directly we see immediately that this term is equal to $1$.

Next, the pole at $s_1=-s_2-s_3-s_4$ gives a residue equal to
\[- \frac{1}
    {\left( {s_2} + {s_3} \right)^2
      \left( { s_2} + { s_4} \right)^2
      \left( { s_3} + { s_4} \right)^2}.\]
Moving $c_2$ to $-\infty$ we pass poles at $s_2=-s_3$ and $s_2 = -s_4$ and the two residues here cancel each other out and we get a contribution of $0$.

For the pole at $s_1 =-s_2$ we have a residue
\[ \frac{      \left( {s_2} - { s_3} - { s_4}
         \right) \left( { s_2} + { s_3} +
        {s_4} \right)e^{s_3 +  s_4}
 }{{ s_2}^2
      \left( { s_2} - { s_3} \right)
      \left( { s_2} + { s_3} \right)
      \left( { s_2} - { s_4} \right)
      \left( { s_2} + { s_4} \right)
      \left( { s_3} + { s_4} \right)^2} .\]
Moving next $c_2$ to $-\infty$ we pass poles at $s_2=0$, $s_2 = -s_3$ and $s_2=-s_4$. The residue here at $s_2=0$ contributes $0$.  The residue at $s_2=-s_3$ is
\[ \frac{s_4(2s_3 +s_4)e^{s_3+s_4}}{2 {s_3}^3(s_3-s_4)(s_3+s_4)^3},\]
which on moving $c_3$ to the left passes  poles at $s_3=0$ and $s_3=-s_4$; on moving $c_4$ to the left the first pole gives a contribution of $-\frac{1}{4}$ and the second gives $0$. In a similar way we find the pole at $s_2=-s_4$ contributes $0$.

Returning to the poles at $s_2 = -s_3$ and $s_2=-s_4$, we will see presently they contribute $0$. Both poles can be handled similarly, so we consider only $s_2=-s_3$. We move $c_2$ to the left and encounter poles at $s_2=-s_3$ and $s_2=-s_4$. The pole at $s_2=-s_3$ gives a residue
\[ \frac{s_4(-2s_3+s_4)e^{-s_3+s_4} }{2{s_3}^3(-s_3+s_4)^3(s_3+s_4)}. \]
We move $c_3$ to the right to make the exponential factor vanish, and encounter only the pole at $s_3=s_4$ which contributes the negative of the residue
\[ \frac{3+6s_4+2{s_4}^2}{16{s_4}^4}\]
which on moving $c_4$ to the left makes a contribution of $0$. In the same fashion the pole $s_2=-s_4$ makes no contribution.

The poles at $s_1=-s_3$ and $s_1=-s_4$ can be handled as above, and one finds that they  contribute $0$. Hence
\begin{equation} \mathcal{C}_4 = \frac{3}{4} .\label{7.31} \end{equation}
\section{A Generalized Result}

To obtain the best numerical results in applications one needs the truncation level $R$ as large as possible, and this  can be optimized by using different truncations $R_i$ in the correlations. Thus, in place of \eqref{6.1} we consider the sum
\begin{equation}\mathcal{S}(\text{\boldmath$j$}, \text{\boldmath$R$}) =\sum_{n=1}^N \Lambda_{R_1}(n+ j_1)\Lambda_{R_2}(n+j_2) \ldots \Lambda_{R_k}(n+j_k),\label{8.1}\end{equation}
and see as before
\begin{equation}\begin{split}
\mathcal{S}(\text{\boldmath$j$},\text{\boldmath$R$})&= N \sum_{\substack{d_i \le R_i, \ 1\le i\le k \\ (d_r,d_s)| j_s-j_r, \, 1\le r<s\le k }}\frac{\prod_{i=1}^k\mu(d_i)\log\frac{R_i}{d_i}}{[d_1, d_2, \ldots , d_k ]} + O(\prod_{i=1}^kR_i) \\ &
=N T_k(\text{\boldmath$j$},\text{\boldmath$R$}) + O(\prod_{i=1}^kR_i). \end{split}\label{8.2} \end{equation}
We now let
\begin{equation} R_i = N^{\theta_i}, \qquad 1\le i\le k, \label{8.3}\end{equation}
and let $\text{\boldmath{$\theta$}}= (\theta_1,\theta_2,\ldots , \theta_k)$.
Then \eqref{6.5} becomes
\begin{equation}  T_k(\text{\boldmath$j$},\text{\boldmath$R$}) = \frac {1}{ (2\pi i)^k}\mathop{\int }_{(c_k)\ }\! \cdots \!\mathop{\int}_{(c_1)\ } F(s_1,\ldots, s_k)\prod_{i=1}^k\frac{N^{\theta_i s_i}}{ {s_i}^2} ds_i ,\label{8.4}\end{equation}
The analysis is now identical to before with the additional constants $\theta_i$ only effecting the constants obtained in the main terms. Suppose as before $\boldsymbol{j} = (j_1,j_2, \cdots j_k)$ has $r$ distinct values which we may take to be the first $r$ components of $\boldsymbol{j}$ with multiplicities $\boldsymbol{a} = (a_1,a_2,\cdots ,a_r)$. Then the truncations $\boldsymbol{\theta} =  (\theta_1,\theta_2,\cdots \theta_k)$ will be partitioned into $r$ sets. We denote by $\theta_i(\nu)$, $1\le \nu \le a_i$ the $\theta$'s associated with $j_i$, and write
\[ \boldsymbol{\theta} = (\boldsymbol{\theta}_1, \boldsymbol{\theta}_2,
\cdots , \boldsymbol{\theta}_r), \quad \boldsymbol{\theta}_i =(\theta_i(1), \theta_i(2), \cdots, \theta_i(a_i)) .\]
Let
\begin{equation} \mathcal{ S}_k(N,\text{\boldmath$j$}, \text{\boldmath $a$},\text{\boldmath{$\theta$}}) =\sum_{n=1}^N \prod_{i=1}^{r}\left(\prod_{\nu=1}^{a_i}\Lambda_{N^{\theta_i(\nu)}}(n+j_i)\right),
\label{8.6}\end{equation}
\begin{theorem} \label{Theorem8} With the above notation, we have for $0<\theta_i \le 1$, and $\log 2||\boldsymbol{j}|| \ll \log N $
\begin{equation} \mathcal{ S}_k(N,\text{\boldmath$j$}, \text{\boldmath $a$},\text{\boldmath{$\theta$}}) =
 \big(\mathcal{ C}_k(\text{\boldmath$a$},\text{\boldmath{$\theta$}})\gs(\text{\boldmath$j$})+o_{k,\text{\boldmath{$\theta$}}}(1)\big)N(\log N)^{k-r} +O(N^{\theta_1+\theta_2 +\cdots +\theta_k}),\label{8.7}\end{equation}
where
\begin{equation} \mathcal{ C}_k(\text{\boldmath$a$},\text{\boldmath{$\theta$}})
= \prod_{i=1}^r \mathcal{C}_{a_i}(\text{\boldmath{$\theta$}}_i ),\label{8.8} \end{equation}
and
\begin{equation} \mathcal{C}_k(\text{\boldmath{$\theta$}})  = \frac {1}{ (2\pi i)^{k}}\mathop{\int}_{(c_k)\ }\!\cdots \!\mathop{\int}_{(c_{1})\ }\Big( \prod_{m \in \mathcal{P}(k)}(\tau_m)^{(-1)^{\# m+1}}\Big) \prod_{i=1}^k\frac{e^{\theta_i s_i}}{ {s_i}^2} ds_i . \label{8.9}
\end{equation}
\end{theorem}

We see that in the special case  $\text{\boldmath{$\theta$}}=(\theta,\theta, \cdots \theta)$ and $R_i = N^\theta$
that
\begin{equation}   \mathcal{ C}_k(\text{\boldmath{$\theta$}})
   = \theta^{k-1} \mathcal{C}_k.\label{8.10} \end{equation}
In the general case we can compute these constants by the same residue calculations used in the last section. The results depend on the relative sizes for the $\theta_i$'s, and we find using the ordering $\theta_1 \ge \theta_2 \ge \cdots \ge \theta_k$ that $\mathcal{C}_1=1$, and
\begin{equation}   \mathcal{ C}_2(\text{\boldmath{$\theta$}})
   = \theta_2 = \min(\theta_1,\theta_2),\label{8.11} \end{equation}
\begin{equation}   \mathcal{ C}_3(\text{\boldmath{$\theta$}})
   = \theta_2\theta_3 - \frac{1}{4}(\theta_2+\theta_3-\theta_1)^2[\theta_2+\theta_3\ge \theta_1],\label{8.12} \end{equation}
and, letting $A_3 = \theta_2+\theta_3-\theta_1$ and $A_4 = \theta_1-\theta_2-\theta_3+\theta_4$,
\begin{equation}  \mathcal{ C}_4(\text{\boldmath{$\theta$}})
   = \theta_2\theta_3\theta_4 - \frac{1}{4}\theta_4{A_3}^2[A_3\ge 0]
 -\frac{1}{32}A_4\left( {A_4}^2+6A_3A_4+4{A_3}^2\right)[A_4\ge 0].\label{8.13}\end{equation}

\end{document}